\documentclass[12pt]{iopart}
\expandafter\let\csname equation*\endcsname\undefined
\expandafter\let\csname endequation*\endcsname\undefined
\usepackage{amsmath}
\usepackage[title]{appendix}
\usepackage{graphicx}
\usepackage{subcaption}
\usepackage{placeins}

\usepackage{amssymb}
\usepackage{dsfont}
\usepackage{iopams}
\usepackage{esint}

\usepackage{siunitx}
\usepackage{algorithm}
\usepackage{algpseudocode}
\usepackage{booktabs}

\usepackage[outdir=./Images/]{epstopdf}

\usepackage[backend=biber,sortcites]{biblatex}
\bibliography{biblio}

\DeclareMathOperator{\diag}{diag}
\DeclareMathOperator{\prox}{prox}
\DeclareMathOperator*{\argmin}{argmin}

\newcommand{\imi}{\mathrm{i}}

\newcommand*{\lvec}[1]{\bi{#1}}
\newcommand*{\lmat}[1]{{#1}}
\newcommand*{\lele}[1]{{#1}}

\begin{document}

\title[Resolving Full-Wave Through-Wall Transmission Effects in Multi-Static SAR]{Resolving Full-Wave Through-Wall Transmission Effects in Multi-Static Synthetic Aperture Radar}

\author{FM Watson$^1$, D Andre$^2$ and WRB Lionheart$^1$}

\address{$^1$ Department of Mathematics, University of Mancheter, Oxford Road, Manchester, M13 9PL, UK}
\address{$^2$ Centre  for Electronic Warfare,  Information and Cyber, Cranfield University,  Defence Academy of the United Kingdom, SN6 8LA, UK}
\ead{francis.watson@manchester.ac.uk}
\vspace{10pt}

\begin{abstract}
Through-wall synthetic aperture radar (SAR) imaging is of significant interest for security purposes, in particular when using multi-static SAR systems consisting of multiple distributed radar transmitters and receivers to improve resolution and the ability to recognise objects.  Yet there is a significant challenge in forming focused, useful images due to multiple scattering effects through walls, whereas standard SAR imaging has an inherent single scattering assumption.  This may be exacerbated with multi-static collections, since different scattering events will be observed from each angle and the data may not coherently combine well in a naive manner.  To overcome this, we propose an image formation method which resolves full-wave effects through an approximately known wall or other arbitrary obstacle, which itself has some unknown ``nuisance'' parameters that are determined as part of the reconstruction to provide well focused images. The method is more flexible and realistic than existing methods which treat a single wall as a flat layered medium, whilst being significantly computationally cheaper than full-wave methods, strongly motivated by practical considerations for through-wall SAR.
\end{abstract}

%
%
%
%
%

\section{Introduction}
Synthetic Aperture Radar (SAR) is a radar collection and image formation technique in which data from many pulses collected along a flightpath are used to \textit{synthesize} a large aperture and form fine-resolution images from large (\textit{stand-off}) distances\cite{jakowatz2012spotlight,cheney2009fundamentals,Cheney09}.  It has widespread use as an important remote sensing technology, partly due to its all-weather, day and night operability (unlike, for example, optical imaging), including for defence and security purposes, measuring biomass, monitoring sea ice, and in monitoring facilities such as ports for commercial reasons. At lower frequencies of around 1 GHz (through UHF and L-band in the IEEE standard), the radar waves readily penetrate many building structures -- else technologies such as mobile phones and WiFi would not work.  This presents the possibility of using SAR for through-wall imaging, sometimes referred to under umbrella of \textit{Remote Intelligence of Building Interiors} (RIBI)\cite{Horne18}, but not without some significant challenge.  The standard data model used as the basis for SAR image formation uses a single scattering assumption -- namely the Born approximation -- but any waves which have interacted with objects of interest inside a building will have scattered at least three times.  This will generally result in artefacts and a hard to interpret image.  

A second challenge is one of data coverage and bandwidth, in order to form a fine enough resolution -- particularly if 3D effects and structures need to be resolved.  At lower frequencies one must fly a longer aperture for the same cross-range resolution (that is, resolution in the same direction as a straight flight path), and frequency bandwidth will inevitably be restricted at lower centre frequencies due to engineering limitations, similarly resulting in a coarser range resolution.  We must also collect data from multiple heights (or elevation angles) in order to gain the 3D data required\cite{Corbett18}, but it may not always be possible to carry out such a complete collection in practice.  One technological solution is to use a \textit{multi-static} radar configuration -- multiple distributed radar transmitters and receivers operating coherently together\cite{santi2014point, Andre21}. A naive approach to combining this multi-static data to form a single image however may simply exacerbate the problem of artefacts and image interpretability: data collected from different source-receiver (or \textit{bi-static}) pairs will have undergone different delays and multiple-scattering through the walls, and so may not coherently combine together naively with a single-scattering model.  

Thus, the need to resolve through-wall full-wave effects and multiple scattering is reinforced.  If we wish to resolve multiple scattering in the inverse problem, then it is reasonable to assume our forward model must include these effects.  Several methods have been developed to form focused images through walls by modelling the through-wall transmission as being through a single- or multi-layered medium\cite{Solimene09,Li10,Zhang14}.  In particular, Solimene \textit{et al}\cite{Solimene09} apply the Kirchhoff approximation in conjunction with an analytic expression for the Green's function through the dielectric layer (wall), and estimate the parameters (permittivity) of the wall as part of a linear reconstruction.  Li \textit{et al}\cite{Li10} similarly estimate wall parameters of a 3-layered medium in order to form a focused back-projection image.  Alternative methods of estimating the properties of the wall may also be employed to improve the results, for example by extracting the off-wall time-delay measurements in mono-\cite{Protiva11} or a bi-static configurations\cite{Elgy21}. The main limitation of a layered medium model is it presupposes some geometrical restrictions on the data collection -- that the radar is observing one outer wall somewhat obliquely. If data is collected from stand-off ranges with airborne radars then we are most likely to find the full structure width within the antenna footprint, and with this reflections from two walls at different angles. If we are to employ multi-static collections then applicability of the layered medium models may be further reduced.

If we do not linearise the inverse problem, for example using some (approximate) Green's function for the medium about which we have linearised, then we must repeatedly solve the forward (wave) model during a non-linear reconstruction process.  A method to do so by solving a series of 1-D inverse scattering problems, which are then combined to form a 2D image, was developed by Klibanov \textit{et al}\cite{Klibanov21}.  As with the layered medium models, this may not apply well to more generalised data collection geometries, since multiple scattering is only resolved along the range direction.  Full-wave inversion techniques, in which the full scattering model (e.g. acoustic or elastic wave equation, Maxwell's equations, etc.) is used in a reconstruction method such as regularised least-squares or Bayesian inversion, have now been widely employed for many inverse scattering problems, in particular for geophysical imaging\cite{Virieux09, Yao20}
and the related ground-penetrating radar (GPR) problem\cite{van2018gpr,klotzsche2019review, Iglesias23,watson2016}.  For through-wall radar, full-wave inversion approaches have been applied in conjunction with level-set techniques both for 2D\cite{Incorvaia19, Incorvaia20} and 3D image formation\cite[ch. 3]{Incorvaia21}. In both of these cases the wall was assumed known (which is reasonable), and somewhat complete data coverage of transmitter/receiver antennas surrounding the building in the near-field was simulated.  The latter may make this challenging to translate to the stand-off SAR imaging problem.  Full-wave inversion also has a significant computational cost due to the need to repeatedly solve wave equations, which may limit its applicability to intelligence gathering applications of SAR in the near-term.

One answer to the computational cost has been the use of reduced order models (ROMs), in which a low-order surrogate model is developed which still encapsulates the degrees of freedom of the full physics (wave) model\cite{lucia2004reduced,Seoane20,Elgy23}.  ROMs have been applied extensively to inverse scattering problems as a means to deal with the nonlinearity due to multiple scattering in various settings by Borcea \textit{et al}\cite{Borcea18, Borcea20, Borcea21lossylayered, Borcea21waves}, as well as more recently to the case of mono-static (i.e. co-located source and receiver) SAR by Druskin \textit{et al}\cite{Druskin22}.  The results provided in this latter work are more akin to a mono-static GPR setup (demonstrating a clear ability to reduce artefacts due to multiple reflections in this context), with the flight path in the imaging plane itself (while the radar will usually be some way above the ground and imaging plane in airborne SAR), close to objects of interest, and with a Gaussian pulse as source.  It is possible there will be some increased computational cost associated with having a fine enough timestep to simulate the wideband frequency-modulated pulses more commonly employed in airborne SAR (or, equivalently, a short enough impulse which can then be deconvolved), as well as in extending the simulation domain to representative ranges.  It may also be necessary to employ a full 3D implementation of the scheme for it to be applicable to representative SAR data.  Notwithstanding the application of interest here being distinctly 3D, in applications of SAR one is most often flying above the image plane, at some elevation angle \(\sigma\) from the scene centre.  Two scatterers in the ground plane at a distance apart \(x\) in the range-direction will appear in the time-domain data \(\tfrac{x\cos(\sigma)}{c}\) \unit{\second} apart, but the multiple reflections will still appear at additional delays of \(\tfrac{x}{c}\) \unit{\second}.  This makes it clear that the problem of dealing with multiple scattering in SAR is a distinctly three-dimensional one -- regardless of whether 2D or volumetric images are being formed.

In this paper, we propose a through-wall imaging method which is able to resolve the full-wave effects through an arbitrary wall (or other obscurant) by modelling the transmission with a boundary element method.  This, combined with a simple reduced order model, is used to form numerical approximations to the Green's functions for an otherwise empty scene, which can be used in a linear reconstruction scheme to image the building interior.  Thus, we can resolve multiple scattering between the wall and objects in the scene, but not between objects in the scene themselves.  By defining the image implicitly as a function of the parameters describing the wall, we obtain a non-linear inverse problem in just a small number of parameters.  This provides us with a scheme which both captures many of the data features observed in real-world relevant data collection scenarios without the above-mentioned limitations of a layered-model approximation. It also has a lower computational cost than fully non-linear (full-wave inversion) methods, more readily allowing application in the short timescales which might be required in building intelligence gathering operations.

The rest of the paper is organised as follows.  First, in Section~\ref{sec: SAR scattering} we discuss the standard data model used in SAR in some detail, to make clear assumptions and how they may cause difficulties, as well as how this model is most commonly used to form SAR images.  In Section~\ref{sec: through wall model} we develop our through-wall scattering model, before discussing its ROM approximation in Section~\ref{sec: ROM}.  Then, in Section~\ref{sec: reconstruction scheme} we present our reconstruction scheme, which we demonstrate in numerical experiments in Section~\ref{sec: numerical results}.  Further details on the boundary element method used for the forward model are included in Appendix~\ref{ap: BEM}, with further details on the reconstruction scheme also being included in Appendix~\ref{ap: reconstruction}.

\section{Single Scattering and Synthetic Aperture Radar}\label{sec: SAR scattering}
\subsection{A standard scattering model for SAR}
While electromagnetic wave propagation is described by Maxwell's equations, the majority of SAR image formation assumes the radar wave propagation is governed by a scalar wave equation\cite{Cheney09,cheney2009fundamentals},
\begin{equation}
	\left(\nabla^2 - \frac{1}{c^2(\mathbf{x})}\frac{\partial^2}{\partial t^2}\right)u(\mathbf{x},t) = j(\mathbf{x},t),
	\label{eq: scalar wave}
\end{equation}
for source \(j\).  This is valid for each component of the electromagnetic waves when travelling through free space, but will not take into account polarization changes upon scattering.  If multi-polar measurements are taken, separate images are generally formed for each polarization channel, which are often colourised and overlaid (or an alternative polarimetric decomposition shown)\cite{Martines21}, which can be used to interpret the type of scattering occurring in the scene.  This is in contrast to a joint reconstruction which determines some anisotropic scattering properties directly, for example as has been carried out for capacitance imaging of liquid crystals\cite{polydorides04}.

Let us decompose the total wavefield \(u\) into the sum of an incident field \(u^{in}\), the field which would be observed in some background medium with speed \(c_0\) but otherwise no objects (possibly a constant), and a scattered field \(u^{sc}\), \(u=u^{in} + u^{sc}\).  Then the incident field satisfies
\begin{equation}
	\left(\nabla^2 - \frac{1}{c_0^2}\frac{\partial^2}{\partial t^2}\right)u^{in}(\mathbf{x},t) = j(\mathbf{x},t),
	\label{eq: incident wave},
\end{equation}
and subtracting (\ref{eq: incident wave}) from (\ref{eq: scalar wave}) we find that the scattered field satisfies
\begin{equation}
	\left(\nabla^2 - \frac{1}{c_0^2}\frac{\partial^2}{\partial t^2}\right)u^{sc}(\mathbf{x},t) = -V(\mathbf{x})\frac{\partial^2 u(\mathbf{x},t)}{\partial t^2},
	\label{eq: scattered wave}
\end{equation}
with reflectivity function \(V(\mathbf{x})=\frac{1}{c_0^2} - \frac{1}{c^2(\mathbf{x})}\).  The solution to (\ref{eq: scattered wave}) can be expressed in terms of the Green's function for the background medium with wavespeed \(c_0\) as
\begin{equation}
	u^{sc}(\mathbf{x},t) = \iint g_0(\mathbf{x}-\mathbf{z}, t-\tau)V(\mathbf{z})u(\mathbf{z},\tau)\, \rmd\tau d\mathbf{z},
	\label{eq: Lippmann time}
\end{equation}
or in the frequency domain for the equivalent Helmholtz problem as
\begin{equation}
	U^{sc}(\mathbf{x},\omega) = -\int G_0(\mathbf{x}-\mathbf{z},\omega)V(\mathbf{z}) \omega^2U(\omega, \mathbf{z})\,\rmd\mathbf{z},	
	 \label{eq: Lippmann freq}
\end{equation}
where \(U^{sc}\), \(U^{in}\), and \(U\) are the Fourier transforms of \(u^{sc}\), \(u^{in}\), and \(u\), respectively.  If \(c_0\) is a constant wavespeed, then the Green's function is given by
\begin{equation}
 (\nabla^2 + k^2)G=-\delta(\mathbf{x}),\qquad	G_0(\mathbf{x},\omega) = \frac{\rme^{\imi k |\mathbf{x}|}}{4\pi|\mathbf{x}|},
 \label{eq: Green's function}
\end{equation}
with wavenumber \(k=\omega/c_0\).

Noting that scattered wavefield appears both in and outside of the integral in (\ref{eq: Lippmann time}-\ref{eq: Lippmann freq}), these do not provide a straightforward formula for the scattered fields which instead vary nonlinearly with \(V\).  As such, SAR (and much radar theory and practice in general) instead makes use of the Born approximation,
\begin{equation}
	U^{sc}\approxeq U_B^{sc} := -\int G_0(\mathbf{x}-\mathbf{z},\omega)V(\mathbf{z}) \omega^2U^{in}(\omega, \mathbf{z})\,\rmd\mathbf{z}.
	\label{eq: Born}
\end{equation}
This is a \emph{single scattering} approximation, and can be justified by expanding (\ref{eq: Lippmann freq}) in a Neumann series and truncating to first order for sufficiently weak scattering.  

\subsection{The SAR data model and image formation}
SAR involves flying a radar along some flight path \(\gamma_T(s)\) during which multiple wide bandwidth pulses are emitted and the echoes received at a receiver radar antenna located at (often co-located) \(\gamma_R(s)\). The parameter of distance along the flight path \(s\) is referred to as slow-time, distinguishing it from the time-of-flight of a pulse \(t\): due to the pulse travelling much faster at the speed of light, one generally assumes everything is stationary during the time of a pulse -- the \textit{start-stop} approximation\cite{Cheney09,cheney2009fundamentals}.   This forms a \textit{synthetic} aperture, providing the two dimensions of information required to form an image: slow time \(s\), or position of antennas; and time of flight at which reflections are received, or frequency \(\omega\).  To form a 3D image, one may either combine the data from multiple data collections at appropriately different apertures (for example different heights), assuming the scene remains stationary, or use a collection multiple transmitters and/or receivers in different locations (a \textit{multi-static} collection).

A simplified model of the incident wavefield in the frequency domain, for a small, point-like antenna, is a multiple of the freespace Green's function, 
\begin{equation}
	U^{in}(\mathbf{z},\omega; s) = P(\omega) \frac{\rme^{\rmi \omega|\mathbf{z}-\gamma_T(s)|/c_0}}{4\pi|\mathbf{z}-\gamma_T(s)|}
	\label{eq: incident}
\end{equation}
for some \(P\in\mathds{C}\) describing the transmitted power and phase at frequency \(\omega\). Thus, the Born approximated scattered field measured at the receiver is given by
\begin{equation}
	U^{sc}_B(s, \omega) = a(\omega)\int \frac{\rme^{\rmi \omega(R_T(s,\mathbf{z}) + R_R(s,\mathbf{z}))/c_0}}{(4\pi)^2R_T R_R}V(\mathbf{z})\,\rmd \mathbf{z},
	\label{eq: Born data}
\end{equation}
where the \(\omega^2\) term has been absorbed into \(a=\omega^2 P\), \(R_{T}(s,\mathbf{z}) = |\mathbf{z} - \gamma_{T}(s)|\), and \(R_{R}(s,\mathbf{z}) = |\mathbf{z} - \gamma_{R}(s)|\).

The received signal will need to be sampled at least at Nyquist frequency. Since the the pulses will be generally have a bandwidth which is small compared to the centre frequency, i.e. a low fractional bandwidth (as an example, the current generation of ICEYE SAR satellites have a centre frequency of 9.65 GHz and bandwidth of 600 MHz\cite{iceyemissions}), the signals are first I/Q demodulated (or \textit{base-banded}) by mixing with a reference signal,
\begin{equation}
	d_{\mathrm{phd}}(s,\omega) = U^{sc}(s,\omega)\rme^{-\rmi\omega\left(R_T(s,\mathbf{x}_{\mathrm{ref}}) + R_R(s,\mathbf{x}_{\mathrm{ref}})\right)/c},
	\label{eq: phdata freespace}
\end{equation}
where the subscript ``phd'' refers to \textit{Phase History Data}, and \(\mathbf{x}_{\mathrm{ref}}\) is some reference location in the scene.

In order to make use of (\ref{eq: Born}) in an image formation method and data model, with discrete pixels, it must be discretised in some way.  In the standard approach to SAR, this amounts to taking \(V\) which is a distribution of point scatterers arranged on some grid on a 2D surface representing the ground -- flat or otherwise -- or in 3D space, which replaces the integral in (\ref{eq: Born}) with a summation.  This effectively assumes that the scene is made of (stationary) isotropic point scatterers.

While there are many different SAR image formation algorithms, most involve applying the adjoint of (\ref{eq: Born}) to the data, or some approximation to it, for example, by making a far-field approximation or by interpolating data-points onto a regular grid in Fourier space to make use of the computationally efficient Fast Fourier Transform (FFT) algorithm\cite{doerry2012polar}. In a mathematical sense these can all be considered a back-projection, but note that ``back-projection'' usually refers to a specific time-domain algorithm in SAR literature\cite{doerry2016bp}.

\section{A forward model for through-wall radar} \label{sec: through wall model}
\subsection{Through-wall transmission}
We noted several assumptions in the standard SAR forward model: governed by a scalar wave equations, only single (weak) scattering occurs, with a scene made of isotropic point scatterers.  Where these assumptions do not hold, artefacts will be observed in the back-projection images based on this forward model.  In some cases artefacts might provide valuable information to an analyst, helping them understand the scattering phenomena which might have occurred and therefore what objects must be present, and so could be considered ``features''. 

In the through-wall imaging application of interest here, the radar waves interacting with objects we would like to image inside a building must have scattered at least 3 times: once in the transmission through the wall; at least once interacting with objects of interest inside the building; then once more in transmission back through the wall. When this data is used to form images, we might see distinct multiple scattering artefacts in the image.  We will almost certainly find that the image does not focus up well -- in particular for multi-static data -- with the different effective path lengths to each locations inside the building preventing coherent combination of the back-projected data.

Consider for now a scene containing a solid object \(D\) with wavespeed \(c_D\) obscuring the radar view, the walls of a building for example, but which is otherwise empty. The scattered wavefield in this case satisfies
\begin{subequations}
\begin{align}
	(\nabla^2 + k_0^2)U^{sc}(\mathbf{x},\omega) &=0, \qquad \mathbf{x}\in \Omega\backslash D, \\
	(\nabla^2 + k_{D}^2)U(\mathbf{x},\omega) &=0 \qquad \mathbf{x}\in D,
\end{align}
\label{eq: Helmholtz inclusion}
\end{subequations}
where \(k_{D}=\omega/c_D\) and \(k_0=\omega/c_0\) are the wavenumbers in free-space and inside the object (respectively), with boundary conditions
\begin{subequations}
\begin{align}
	\gamma_0^+(U^{sc}(\mathbf{x},\omega) + U^{in}(\mathbf{x},\omega))= \gamma_0^- U(\mathbf{x},\omega),  &\qquad \mathbf{x}\in\partial D \\
	\frac{1}{c_0}\gamma_1^+(U^{sc}(\mathbf{x},\omega) + U^{in}(\mathbf{x},\omega)) = \frac{1}{c_D}\gamma_1^- U(\mathbf{x},\omega),&  \qquad \mathbf{x}\in\partial D,\\
\end{align}
\end{subequations}
as well as the Sommerfeld radiation condition. Here, \(\gamma_0^+\) and \(\gamma_0^-\) are the Dirichlet trace operators -- i.e. the limit of values towards the interface \(\partial D\) from the exterior and interior of the object (respectively), and \(\gamma_1^+\) and \(\gamma_1^-\) the Neumann trace operators -- i.e. the limits of the normal derivatives of the wavefield at the interface.  The solution to this equation can be written using Green's representation theorem as a combination of single and double layer potentials\cite{Betcke17}
\begin{equation}
	U(\mathbf{x},\omega) = \int_{\partial D} G_0(\mathbf{x},\mathbf{y},\omega)\psi(\mathbf{y}) \rmd\mathbf{y} - \int_{\partial D} \frac{\partial G_0(\mathbf{x},\mathbf{y},\omega)}{\partial \mathbf{n}}\phi(\mathbf{y}) \rmd\mathbf{y},
	\label{eq: representation formula}
\end{equation}
where \(G_0(\mathbf{x},\mathbf{y},\omega)\equiv G_0(\mathbf{x}-\mathbf{y},\omega)\), and with jumps of the solution across the interface
\begin{align}
	\psi &= \gamma_1^- u - \gamma_1^+ u, \\
	\phi &= \gamma_0^- u - \gamma_0^+ u.
\end{align}

The scattering problem (\ref{eq: representation formula}), being a surface integral, is suitable for solving with the Boundary Element Method (BEM), for which we use the BEMPP-CL library\cite{betcke2021bempp}.  In scattering problems consisting of one or a small number of discrete solid objects (rather than a continuum of varying wavespeed), depending on their size and number, boundary element methods may provide a computational benefit over e.g. a finite element method, due to the need to only discretise and solve equations only on the surface(s) rather than a whole volume. 

The derivation of the specific boundary element formulation of (\ref{eq: reconstruction problem full}) which we use is included in Appendix~\ref{ap: BEM}.

\subsection{The data model}
Now, returning to our otherwise empty scene containing the single obstacle -- a wall, say.  If we take for the incident wavefield \(U^{in}\) in (\ref{eq: combined bie}) the freespace Green's function \(G_0(\mathbf{x},\mathbf{y},\omega)\) (\ref{eq: Green's function}), i.e. with source term \(J(\mathbf{x})=\delta(\mathbf{x}-\mathbf{y})\) with \(\mathbf{y}\in\mathds{R}^3\backslash D\), then the total field 
\begin{equation}
	G'(\mathbf{x},\mathbf{y},\omega):=U(\mathbf{x}) = G_0(\mathbf{x},\mathbf{y},\omega)+U_{D}^{sc}(\mathbf{x},\omega)
\end{equation} 
will be the Green's function \(G'\) for the domain with the obscurant, where we have denoted \(U_D^{sc}\) the wavefield scattered by \(D\).

Let us parameterise the obscuring object \(D\) by some vector \(\lvec{m}\) describing its unknown properties, which may be a list (for example) of electromagnetic properties (permittivity, conductivity) and/or some physical dimensions.  We may thus write \(U_D^{sc}(\mathbf{x},\mathbf{y},\omega):=U_D^{sc}(\lvec{x},\lvec{y},\omega;\lvec{m})\), and \(G'(\mathbf{x},\mathbf{y},\omega):=G'(\mathbf{x},\mathbf{y},\omega;\lvec{m})\).  With this parameterisation in mind, our model for through-wall radar data is obtained by substituting \(G'\) for the free-space Green's function in (\ref{eq: Born}) and (\ref{eq: incident}) (where in the latter it is written explicitly), plus the wavefield scattered directly by \(D\) measured at the antenna.  This results in a Born approximation with corrected Green's function
\begin{equation}
	\lvec{d}_{\mathrm{phd}} = \mathcal{F}(V;\lvec{m}):=\mathcal{F}_0(\lvec{m}) + \mathcal{F}_1(V;\lvec{m}), 
	\label{eq: through wall model}
\end{equation}
which has components
\begin{subequations}
	\begin{align}
		[\mathcal{F}_0(\lvec{m})]_i:= &a(\omega_i)\rme^{-\rmi \omega R_{0,i}/c}U_D^{sc}(\mathbf{x}_{R,i}; \mathbf{x}_{T,i}, \omega_i,\lvec{m}) ,  \\
		[\mathcal{F}_1(V;\lvec{m})]_i:= &a(\omega_i)\rme^{-\rmi \omega R_{0,i}/c}\cdot\nonumber\\ 
		&\int G'(\mathbf{x},\mathbf{y}_{T,i}, \omega; \lvec{m})G'(\mathbf{y}_{R,i} \mathbf{x},\omega; \lvec{m})V(\mathbf{x})\rmd\mathbf{x},
	\end{align}
	\label{eq: through wall model components}
\end{subequations}
where \([\mathcal{F}_{\cdot}]_i\) is the contribution to the i\textsuperscript{th} sample of the  measured data (ordered into a vector), with \(\omega_i\) the i\textsuperscript{th} measured angular frequency, \(\mathbf{x}_{R,i}\) and \(\mathbf{x}_{T,i}\) the i\textsuperscript{th} receiver and transmitter locations, and \(R_{0,i}:=R_T(s_i, \mathbf{x}_{\mathrm{ref}})+R_R(s_i, \mathbf{x}_{\mathrm{ref}})\) for some chosen reference point in the scene \(\mathbf{x}_{\mathrm{ref}}\). As previously, \(V(\mathbf{x})\) is the reflectivity describing objects in the scene \emph{excluding} the obscurant (wall).  Note that \(\mathcal{F}_0\) only contains the scattered field, since the incident field will not be recorded due to antenna beam patterns and time-domain gating of measured signal.

Before discussing numerical evaluation of (\ref{eq: through wall model}), we first re-emphasise what features of the data are modelled -- and therefore we expect to be able to resolve well -- compared to the standard Born approximated forward model.  (\ref{eq: through wall model}) includes (multiple) scattering  of the incoming waves through and off of a (known) wall to the target scene, single scattering off of other targets in the scene, before being (multiply) scattered back through the wall to the receiver.  It also includes  multiple scattering off of the known wall itself without other interactions in the scene.  It does not include multiple scattering between other objects in the scene, which are modelled using the Born approximation.  Other standard assumptions also remain -- such as scalar waves and a point scatterer approximation for objects \emph{in} the scene (not the wall).

We should therefore expect use of (\ref{eq: through wall model} - \ref{eq: through wall model components}) in a reconstruction algorithm to help resolve through-wall defocusing and multiple scattering effects (including off of side walls and corners).  This might, for example, resolve ``ghost'' reflections of objects inside buildings where these are caused by multiple reflections off of an external wall (e.g. side wall), but would not resolve multiple reflections between other objects in the scene. 

Note that including the \(\mathcal{F}_0\) term in a reconstruction scheme would result in the wall not being seen in the reconstructed image (if we were able to perfectly match the reflected waves), so it could be desirable to neglect this term.

\FloatBarrier
\section{Reduced order modelling} \label{sec: ROM}
The model (\ref{eq: through wall model components}) is to be used within an iterative reconstruction for discretised reflectivity \(V:=\sum_i v_i\delta(\tilde{\mathbf{x}}_i)\), \(\lvec{v} = [v_1,\ldots, v_n]^T\), \(v_i = V(\tilde{\mathbf{x}}_i)\).  Clearly re-solving the boundary integral equations at each iteration would present such a computational cost that we may as well have carried out a full-wave reconstruction, but we can pre-calculate the Green's functions \(G'\) for a given \(\lvec{m}\).  

If \(\lvec{m}\) is known, then we could pre-compute numerical Green's functions \(G'\) to be stored and reused. If instead \(\lvec{m}\) is unknown -- as will generally be the case -- we can use a Reduced Order Model (ROM). This involves pre-computing and storing only a number of Green's functions for a predetermined set of parameters (the \textit{offline} or training stage), which can later be used to accurately approximate in real-time the wavefield solutions (the \textit{online} stage).  As well as allowing for a real-time evaluation of the numerical Green's functions, the total number of numerical simulations may also be greatly reduced.

For the ROM, following Seoane \textit{et at}\cite{Seoane20}, we use a Proper Orthogonal Decomposition with Interpolation (PODI) of the simulated \(U_D^{sc}(\tilde{\mathbf{x}}; \lvec{m}_i)\) for a predetermined set of observation parameters \(\{\lvec{m}_i\}\).

Denote by \(\lvec{u}^{sc}(\lvec{m}):= [U^{sc}_{D}(\tilde{\mathbf{x}}_0; \lvec{m}), \ldots, U^{sc}_{D}(\tilde{\mathbf{x}}_n; \lvec{m})]\) the simulated scattered wavefield at locations in the scene to be imaged \(\tilde{\mathbf{x}}_j\), \(j=1,\ldots,N\), for the set of observation parameters \(\lvec{m}\).  We wish to approximate these scattered fields with some set of basis vectors \(\lvec{\Phi}_k\)
\begin{equation}
	\lvec{u}^{sc}(\lvec{m}) \approx \sum_{k=1}^M a_k \lvec{\Phi}_k,
\end{equation}
with \(M\ll N\).  Having carried out prior observations \(\lvec{u}^{sc}_{D,i}:=\lvec{u}_D^{sc}(\lvec{m}_i)\), we assemble an observation matrix
\begin{equation}
	\lmat{D} = [\lvec{u}^{sc}_{D,1}, \ldots, \lvec{u}^{sc}_{D,N_{\mathrm{obs}}}] \in \mathds{C}^{N \times N_{obs}},
\end{equation}
and take its truncated singular value decomposition 
\begin{equation}
	\lmat{D} \approx \tilde{\lmat{D}} := \lmat{H}^K\lmat{\Sigma}^K(\lmat{G}^K)^* = \sum_{k=1}^K \lvec{h}_k\sigma_k\lvec{g}_k^*,
\end{equation}
for \(K\) corresponding to a truncation in singular values \(\sigma_k\geq \sigma_1\alpha\), for some threshold \(0<\alpha\ll 1\), \(\lmat{H}\in\mathds{C}^{N\times K}\), \(\lmat{\Sigma}=\diag_0([\sigma_1,\ldots,\sigma_K]^T)\in\mathds{R}^{K\times K}\), and \(\lmat{G}^K\in\mathds{C}^{N_{obs}\times K}\). Noting that the observations can be approximately recovered as
\begin{equation}
	\lvec{u}^{sc}_i\approx H^K\sigma^K\lvec{g}_i^*,
	\label{eq: POD recovery}
\end{equation}
the PODI scheme is to approximate \(\lvec{u}^{sc}(\tilde{\lvec{m}})\) for a given \(\tilde{\lvec{m}}\) by
\begin{equation}
	\lvec{u}^{sc}(\tilde{\lvec{m}})\approx \lvec{u}^{sc}_{P}(\tilde{\lvec{m}}):= H^M\Sigma^M \lvec{g}_{P}(\tilde{\lvec{m}})^*,
	\label{eq: POD approximation}
\end{equation}
where \(\lvec{g}_{P}(\tilde{\lvec{m}})\) is obtained via cubic interpolation of the \(N_{obs}\) truncated right singular vectors \(\lvec{g}_k\in\mathds{C}^K\), which correspond to the observation interpolation points \(\lvec{m}_k\).

Writing the discretised forward model (\ref{eq: through wall model}) in terms of the POD, we then have
\begin{align}
	\left[\mathcal{F}_0(\lvec{m})\right]_i &= a(\omega_i)\rme^{-\rmi\omega_i R_{0,i}/c} \left[\lvec{u}^{sc}_{P_0}(\lvec{m})\right]_i  \nonumber\\
	&=  a(\omega_i)\rme^{-\rmi\omega_i (R_{T,i} + R_{R,i})/c} \left[H_{P_0}^M\Sigma_{P_0}^M \lvec{g}_{P_0}(\lvec{m}, \omega_i, s_i)\right]_{R,i},
	\label{eq: POD F0}
\end{align}
with subscript \(P_0\) referring to a POD formed of observations associated with \(\mathcal{F}_0\) of \(U_D^{sc}\) at the receiver locations, subscript \(R,i\) denoting the i\textsuperscript{th} receiver index (for a multi-static configuration).  Similarly, \(\mathcal{F}_1\) is given by
\begin{align}
	\left[\mathcal{F}_1(\lvec{m},\lvec{v})\right]_i =& a(\omega_i)\rme^{-\rmi\omega_i R_{0,i}/c}\sum_j\Bigr\{&& \left(G_0(\lvec{x}_j, \lvec{y}_{T,i}) + U^{sc}_{P_1}(\lvec{x}_j, \lvec{y}_{T,i}) \right) \cdot& \nonumber\\
	& &&\left(G_0(\lvec{x}_j, \lvec{y}_{T,i}) + U^{sc}_{P_1}(\lvec{x}_j, \lvec{y}_{T,i}) \right)v_j \Bigr\} \\
	=& a(\omega_i)\rme^{-\rmi\omega_i R_{0,i}/c}\sum_j\Bigr\{&& \left(G_0(\lvec{x}_j, \lvec{y}_{T,i}) + \left[H_{P_1}^M\Sigma_{P_1}^M \lvec{g}_{P_1}(\lvec{m}, \omega_i, s_i)\right]_{T,ij} \right) \cdot& \nonumber\\
	& && \left(G_0(\lvec{x}_j, \lvec{y}_{T,i}) + \left[H_{P_1}^M\Sigma_{P_1}^M \lvec{g}_{P_1}(\lvec{m}, \omega_i, s_i)\right]_{R,ij} \right)v_j \Bigr\}.
	\label{eq: POD F1}
\end{align}
As with (\ref{eq: POD F0}), subscript \(T,i\) denotes the i\textsuperscript{th} transmitter index (in a multi-static configuration), and subscript \(P_1\) refers to a POD formed with observations associated with \(\mathcal{F}_1\) being \(U_D^{sc}\) at points in the scene \(\mathbf{x}_j\).  This model can be written more succinctly as
\begin{equation}
	\mathcal{F}[\lvec{m},\lvec{v}] := \mathcal{F}_0(\lvec{m}) + \lmat{A}(\lvec{m})\lvec{v},
\end{equation}
where the elements of the matrix \(\lmat{A}\in\mathds{C}^{(n_{\omega}n_{s})\times n}\) are defined by (\ref{eq: POD F1}).

Recall that the PODI scheme is carried out in two steps.  First, the \textit{offline stage} involves prior simulation of \(\lvec{u}^{sc}_{D,k}\) and calculation of \(\lmat{H}^K\), \(\lmat{\Sigma}^K\), and \(\lmat{G}^K\) for the set of interpolation points \(\{\lvec{m}_i\}_{i=1}^{N_{obs}}\).  The offline stage is carried out prior to reconstruction, and the terms stored.  Then, the \textit{online stage} involves evaluating (\ref{eq: POD approximation}) during the iterative reconstruction process, instead of evaluating (\ref{eq: bempp system}) with (\ref{eq: through wall model}) directly.

We note that more advanced projection-based POD methods (or PODP) are also available, which are applied directly to the finite element system, reducing the size of this system to be solved\cite{Seoane20}. These have several advantages compared to PODI, including forcing the solutions to obey the equation of the problem, often yielding greater accuracy and robustness or allow fewer (adaptively chosen) observation points. Our choice is made primarily for implementational simplicity.  

The choice of a PODI scheme being applied directly to the wavefield in the imaging domain also reduces the computational cost of the online stage, since evaluation of the representation formula (\ref{eq: representation formula}) has already been carried out. ROM schemes applied to the BEM solution vector (including PODP) would require evaluation of this representation formula during the online stage, increasing the computational cost during a reconstruction.
This does come at a cost of memory: for a scene with many pixels the observations are much larger than the BEM solution (which scales with size of wall/object, not the scene size), and we must apply an SVD to the whole observation matrix \(\lmat{D}\).

\FloatBarrier
\section{Through-wall reconstruction scheme} \label{sec: reconstruction scheme}
\subsection{Variable projection}
The reconstruction problem in both the model parameters \(\lvec{m}\) and scene \(\lvec{v} := [v_1,\ldots, v_n]^T\), \(v_i = V(\tilde{\mathbf{x}}_i)\), can be formed as the variational problem
\begin{equation}
\begin{aligned}
	\tilde{\lvec{m}}, \tilde{\lvec{v}} =& \argmin_{\lvec{m},\lvec{v}} \mathcal{J}(\lvec{m},\lvec{v}), \\
	\mathcal{J}(\lvec{m},\lvec{v}):=& \frac{1}{2}\|\mathcal{F}(\lvec{m},\lvec{v}) - \lvec{d}\|_2^2 + \lambda_m\mathcal{M}(\lvec{m}) + \lambda_v \mathcal{R}(\lvec{v}),
\end{aligned}
\label{eq: reconstruction problem full}
\end{equation}
with \(\mathcal{M}\) and \(\mathcal{R}\) some regularisation terms enforcing prior knowledge and preventing over-solving noise (which may also include hard constraints), and \(\lambda_m,\lambda_v>0\).  For our purpose, the wall parameters \(\lvec{m}\) may be considered \textit{nuisance parameters}\cite{aravkin2012estimating}, since the presence of a wall is self evident. (\ref{eq: reconstruction problem full}) is well-suited to \textit{reduced modelling} or \textit{variable projection},
\begin{subequations}
\begin{align}
	\tilde{\lvec{m}} &= \argmin_{\lvec{m}} \bar{\mathcal{J}}(\mathbf{m}) \label{eq: projected base}\\
	\bar{\mathcal{J}}(\lvec{m})&:= \mathcal{J}(\lvec{m};\bar{\lvec{v}}(\lvec{m})) \label{eq: projected m}\\
	\bar{\lvec{v}}(\lvec{m}) &:= \argmin_{\lvec{v}} \mathcal{J}(\lvec{m},\lvec{v}). \label{eq: projected v}
\end{align}
\label{eq: reconstruction problem projected}
\end{subequations}
The solutions to (\ref{eq: reconstruction problem full}) and (\ref{eq: reconstruction problem projected}) clearly coincide, simply by back-substitution of (\ref{eq: projected m}) and (\ref{eq: projected v}) into (\ref{eq: projected base}).  However, variable projection here has the benefit of having separated a simpler linear reconstruction problem (\ref{eq: projected v}) from the full non-linear one (\ref{eq: reconstruction problem full}). This approach enables a wide range of black-box optimisation schemes to be applied to the optimisation sub-problems, including cases where the regularisation term \(\mathcal{M}\) is non-smooth\cite{van2016non,van2021variable}. Details of the derivatives of the objective function (\ref{eq: reconstruction problem projected}) used in the reconstruction are provided in Appendix~\ref{ap: reconstruction derivs}.

Note that, unlike in the afore referenced work \cite{aravkin2012estimating}, in (\ref{eq: reconstruction problem projected}) it is the nuisance parameters \(\lvec{m}\) which the data misfit function varies non-linearly with, rather than the parameters of interest \(\lvec{v}\).  

\subsection{Numerical scheme}
In reconstructing the complex-valued reflectivities \(\lvec{v}\in\mathds{C}^n\), we apply regularisation only to their magnitudes which is the  quantity of interest, since the phase of the image is driven by positional error of the voxel versus true scatterer locations (or them not being well represented by point scatterers) and shall vary freely in \([0,2\pi)\).  A suitable and convenient choice is Total Variation (TV), since the interior of a building will presumably consist of some distinct solid objects and be otherwise empty -- which results in the inner optimisation problem for \(\lvec{v}\)
\begin{equation}
	\bar{\lvec{v}}(\lvec{m}) := \argmin_{\lvec{v}}\|\lmat{A}(\lvec{m})\lvec{v}-\left(\lvec{d}-\mathcal{F}_0(\lvec{v})\right)\|_2^2 + TV(|\lvec{v}|),
	\label{eq: implicit v}
\end{equation}
where \(|\cdot|\) is understood to mean the element-wise absolute value \(|\lvec{v}| = [|v_1|,\ldots,|v_n|]^T\), and \(TV(\cdot)\) is the (isotropic) Total Variation semi-norm.  Details of TV and numerical solution to (\ref{eq: implicit v}) are provided in Appendix~\ref{ap: algorithm reflectivity}

Given implicitly defined \(\bar{\lvec{v}}(\lvec{m})\), the nonlinear problem in nuisance parameters (\ref{eq: projected m}) can be solved with a suitable non-linear optimisation scheme, for which we have chosen the BFGS scheme. Again, further implementational details are provided in Appendix~\ref{ap: algorithm outer}.

\FloatBarrier
\section{Numerical results} \label{sec: numerical results}
Here we test the method with simulated data in three numerical experiments: a reconstruction with a perfectly known wall; a reconstruction with a wall of unknown permittivity, but otherwise known geometry; and a reconstruction for a wall of unknown permittivity, and approximately known position and thickness. For the latter, the data is also simulated fully with a full-wave model. Thus, we progressively reduce the \textit{inverse crime} committed to demonstrate the potential of this method.

In each of these numerical experiments, the obstacle comprises two joining outer walls at right angles to one-another as a surrogate for the facing corner of a building (referred to as a ``corner wall'' for brevity), which is designed to mimic our previous experimental setups\cite{Andre21, andre22}.  

\begin{figure}
	\centering
	\includegraphics[width=0.6\textwidth]{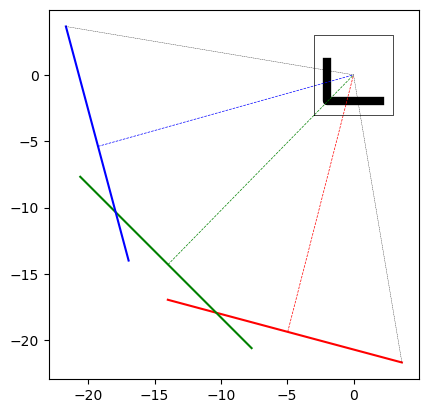}
	\caption{Experimental multi-static collection geometry used for simulations relative to the wall location and image extent.  The flight paths of the transmit and first receiver are shown in red (solid), of the second receiver in green, and of the third receiver in blue. The angles between the lines bisecting the flightpath to the scene centre (dotted lines) are given in table~\ref{tab: radar params}, referred to as the ``bi-static angles''.  Also highlighted by black-dotted lines are the lines of sight from the widest extent of antenna positions to the scene centre.}
	\label{fig: geometry}
\end{figure}

The geometry of the multi-static collection relative to the scene and corner wall are shown in \figurename~\ref{fig: geometry}, in which we ensure there is no line-of-sight from objects in the scene to the antenna positions.  The radar collection parameters are also provided in table~\ref{tab: radar params}.  The angle between the lines bisecting the centre of transmit and receiver apertures to the scene centre are referred to in this table as the ``bi-static angles'', with the multi-static datasets consisting of all 3 transmitter-receiver pairs, including the co-located pair.

Each reconstruction uses data simulated using a 3D forward model, but a 2D dataset is collected (i.e. there is no vertical aperture), and 2D images of the scene are formed using the 3D forward model derived in Section~\ref{sec: through wall model}.  This model includes out-of-plane through-wall transmission and reflection effects.  The same method may be directly applied to 3D reconstructions to also resolve objects in height, simply at increased computational cost, and at the requirement of a vertical aperture of data also being collected.

The centre frequency listed in table~\ref{tab: radar params} may be considered slightly low for a practical through-wall SAR system, but it is of the right order to be representative, and allows us to keep the number of boundary elements low for more rapid testing.  Similarly, the simulated range of \qty{20}{\meter} can easily be increased or decreased as required without any change to computational cost and (given the same extent of aperture angle) little to no affect on the reconstruction.

\begin{table}[h]
	\centering
	{\footnotesize
\begin{tabular}{l |c}\toprule
		Parameter	& Simulation value  \\\midrule
			Centre frequency (\unit{\mega\hertz})& 349.9 \\
			Bandwidth (\unit{\mega\hertz}) & 299.8 \\
			Azimuth aperture (\unit{\radian})& 0.86 \\\midrule
			Transmitter azimuth from \(\hat{\mathbf{x}}\) (\unit{\radian}) & \(-7\pi/12\) \\ 
			Bi-static angles (\unit{\radian})& \(0\), \(-\pi/6\), \(-\pi/3\)\\
			Range to scene (\unit{\meter}) & 20 \\\midrule
			Horizontal pixel spacing (\unit{\meter}) & \num{0.05} \\
			\parbox[m]{5cm}{Effective mono-static \\cross-range resolution (\unit{\meter})} & \num{0.5}\\
			\parbox[t]{5cm}{Effective mono-static \\range resolution (\unit{\meter})} & \num{0.5} \\
			\bottomrule
\end{tabular}
		\caption{Radar data collection and image formation parameters used in simulations.}
		\label{tab: radar params}}
\end{table}

\subsection{Reconstruction with a known obstacle}\label{sec: 2d crime}
To ensure the method is capable of focusing through-wall effects, we first compare the reconstruction results using the through-wall data model with the exact permittivity and geometry used in the data simulation to that of the standard SAR model.  

Using the model (\ref{eq: POD F1}), data is simulated for three point scatterers of reflectivity \num{1.0} beyond the corner wall of height \qty{3}{\meter}, wall lengths \qtylist{3;4}{\meter}, width \qty{0.3}{\meter} and relative permittivity \num{3.}.  These point scatterers do not lie exactly on pixel locations of the reconstruction domain, which themselves are upsampled by a factor of 10 versus the \qty{0.5}{\meter} range- and cross-range resolution of a single (mono-static) channel of the multi-static data collection.  For a more transparent comparison between the two operators, the energy directly reflected from the wall given by \(\mathcal{F}_0\) is not included in the data, i.e.
\begin{equation}
	\lvec{d}_{\mathrm{phd}} = \lmat{A}(\lvec{m})\lvec{v}.
\end{equation}

\begin{figure}[ht]
	\centering
	\begin{subfigure}[b]{0.45\textwidth}
		\centering
		\includegraphics[width=\columnwidth]{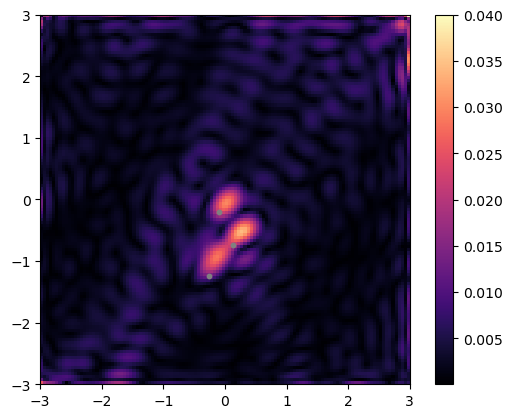}
		\caption{Standard SAR model}
	\end{subfigure}~
	\begin{subfigure}[b]{0.45\textwidth}
		\centering
		\includegraphics[width=\columnwidth]{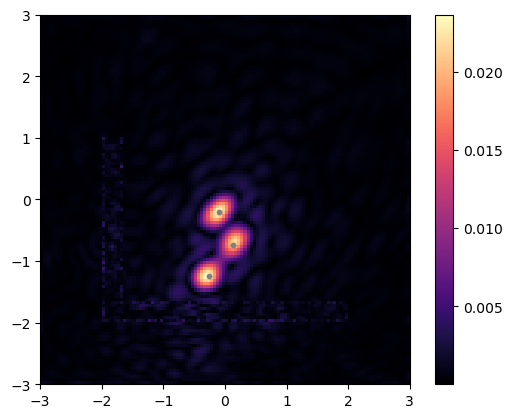}
		\caption{Through-wall SAR}
	\end{subfigure}~
	\caption{Reconstruction results of through-wall data for three point scatterers, using (a) the standard SAR model and (b) the exact through-wall model.  The true location of scatterers are overlaid as grey dots, and the  corner of the wall is at \((-2,-2)^T\).}
	\label{fig: 2d recon inverse crime}
\end{figure}

For both reconstructions, we allow FISTA to run for 500 iterations, and set the regularisation parameter \(\lambda_v\) to zero (i.e. \(\mathcal{R}\equiv 0\)) to test only the ability to fit the noiseless data.  The reconstruction results are shown in \figurename~\ref{fig: 2d recon inverse crime}, with the exact location of the point scatterers overlaid in red.  By comparing these images, the effect of through-wall delay using the standard SAR model is clear, with the targets shifted down-range.  They also appear less well focused and separated, with significantly more background clutter.  In contrast, using the (exact) through-wall SAR model, the location of objects is corrected, it is clearer that there are three distinct objects, and there is a clearer background (i.e. reduced ``clutter'').

\begin{figure}[ht]
	\centering
		\includegraphics[width=0.5\columnwidth]{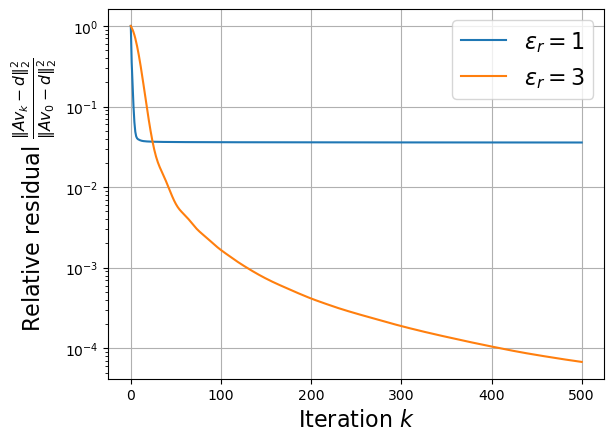}
	\caption{Convergence of the reconstructions using FISTA, for the standard SAR model \(\epsilon_r=1\) shown in blue, and the exact through-wall model \(\epsilon_r=3\) in orange.}
	\label{fig: 2d recon inverse crime convergence}
\end{figure}

\figurename~\ref{fig: 2d recon inverse crime convergence} shows the convergence of FISTA, from which we can see the standard SAR reconstruction quickly reaches a minima in which the data fit is not good, whereas the through-wall SAR model continues to converge.  The slower rate of convergence for the through-wall model is due to the standard choice of step size in FISTA being the reciprocal of the Lipschitz constant of the data misfit, but the operator norm of \(\lmat{A}(\epsilon_r)\) was observed to increase with \(\epsilon_r\).

Note that in both cases, we would not expect either of these reconstructions to reproduce the correct reflectivity value \num{1.0}.  Since the images are up-sampled versus the resolution we could expect from the data, this distributes a single scatterer (and its scattering power) across multiple pixels.

\begin{figure}[ht]
	\centering
	\begin{subfigure}{0.3\textwidth}
		\centering
		\includegraphics[width=0.9\textwidth]{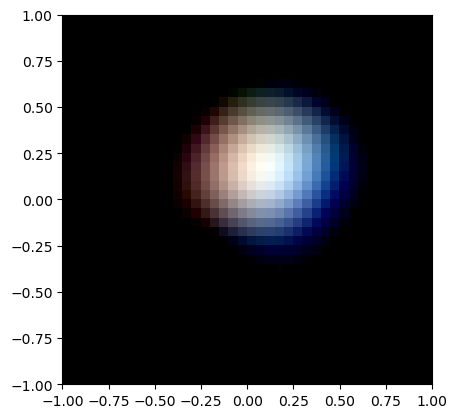}
		\caption{Standard SAR model}
		\label{fig: colourised std sar}
	\end{subfigure}~~
	\begin{subfigure}{0.3\textwidth}
		\includegraphics[width=0.9\textwidth]{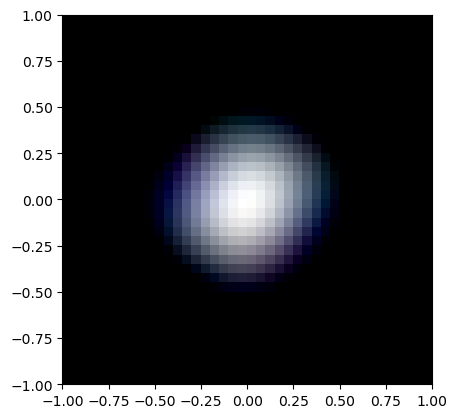}
		\caption{Through-wall SAR}
		\label{fig: colourised tw sar}
	\end{subfigure}~~
	\begin{subfigure}{0.25\textwidth}
		\centering
		\includegraphics[width=0.9\textwidth]{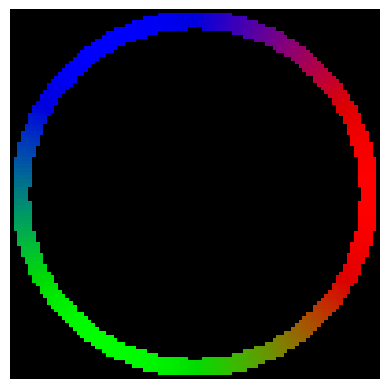}
		\caption{Colour wheel}
	\end{subfigure}
	\caption{Overlaid magnitude of back-projection images for a single point target at the origin behind a corner wall. A separate image is formed from data of each bi-static pair in the multi-static collection. These are colourised and overlaid such that pixels taking the same value in each image which appear grey-scale, or otherwise tinted towards red, green or blue if the image is brighter in the corresponding image. Red, green and blue image colour channels also correspond to the antenna flight paths shown in \figurename~\ref{fig: geometry}.}
\end{figure}

To understand where the difference performance arises, it is useful to consider the case of a single point scatterer.  For this case, individual back-projection images are formed for each transmit/receiver pair of data in the multi-static collection. These 3 images are then colourised and overlaid such that pixels taking the same value in each image which appear grey-scale, and will otherwise be tinted towards red, green or blue if the value is greater in the first, second or third channel's image.  The results are shown in \figurename~\ref{fig: colourised std sar} for images formed using the standard SAR model, and in \figurename\ref{fig: colourised tw sar} using the through-wall model with the correct permittivity.  

The increased tinting from red to blue observed in \figurename~\ref{fig: colourised std sar} is because the peaks of each channel's image are not co-located, resulting also in phase errors between the images due to different effective path lengths through the wall to each pixel.  This means that when all the back-projected data is summed it will not properly coherently combine.  In contrast, we see less colourisation in \figurename~\ref{fig: colourised tw sar} due only to different sidelobe patterns, as the peaks of each image are properly aligned since the effective path lengths corrected. This suggests the multi-static data contains information which may be used to help resolve an unknown wall in order to align each data channel.

\begin{figure}[ht]
	\centering
	\begin{subfigure}{0.9\textwidth}
		\centering
		\includegraphics[width=0.8\textwidth]{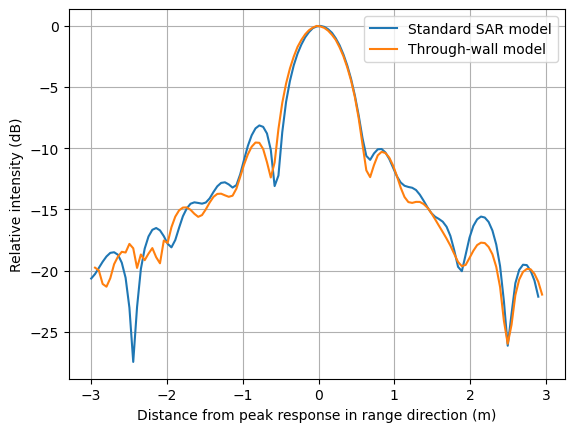}
		\caption{Sidelobe pattern of backprojection images in the range direction.}
		\label{fig: range sidelobe}
	\end{subfigure}~\\~\\~\\
	\begin{subfigure}{0.9\textwidth}
		\centering
		\includegraphics[width=0.8\textwidth]{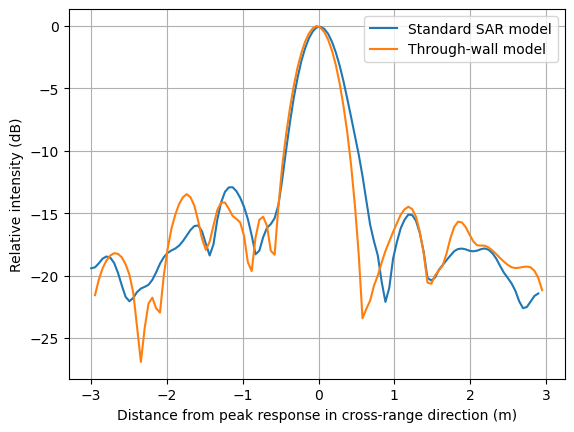}
		\caption{Sidelobe pattern of backprojection images in the cross-range direction.}
		\label{fig: xrange sidelobe}
	\end{subfigure}
	\caption{Sidelobe pattern of backprojection images scaled to peak amplitude of the main lobe, shown in log scale. Blue shows the standard SAR model, and orange the through-wall SAR model.}
	\label{fig: sidelobes}
\end{figure}

In \figurename~\ref{fig: sidelobes}, we also plot cross-sections through the back-projection images for the whole dataset to show the sidelobe patterns of a single point scatterer.  \figurename~\ref{fig: range sidelobe} shows nominally equal resolution for both the standard SAR and through-wall models in the range direction for this geometry.  However, some sidelobe performance is gained using the through-wall model, with approximately \SI{1.5}{\decibel} improvement in peak sidelobe performance.  \figurename~\ref{fig: xrange sidelobe} demonstrates some improved resolution in the cross-range direction achieved by using through-wall model.  Whilst there is still a gain of approximately \SI{1.0}{\decibel} peak sidelobe performance, there is less gained overall.  In either case, no apodization has been performed to control sidelobe levels.

Note a single point scatterer is an unrealistic ``best case'' when using the standard SAR model, as it may only appear slightly displaced and defocused, without interference of responses from other objects overlaying due to multiple angles of observation.

\subsection{Reconstruction with an unknown permittivity}\label{sec: 2d partial crime}
Now we consider the case where the permittivity of the wall is unknown.  The data simulation and experimental setup is the same as for Section~\ref{sec: 2d crime}, except for the true wall permittivity and thickness being reduced to \num{2.07} and \qty{0.2}{\meter} respectively. \qty{20}{\percent} Gaussian white noise is added to the data relative to its norm. 

We set the regularisation parameter for TV to \(\lambda=10^{-3}\|\lmat{A}\|_2\), for \(\lmat{A}\) corresponding to \(\epsilon_r=3.0\), with this norm being calculated using the power method. This choice of \(\lambda\) being comparatively small to \(\|\lmat{A}\|_2\) is so as not to drive the reconstruction.  As before, we use only the operator \(\mathcal{F}_1\) in simulation and reconstruction (as if the first reflections from the wall had been filtered out of the received data), so we can observe how well \(\lvec{m}\) is resolved using only information from objects beyond the wall. 

\begin{figure}[ht]
	\centering
	\begin{subfigure}[b]{0.45\textwidth}
		\centering
		\includegraphics[width=\columnwidth]{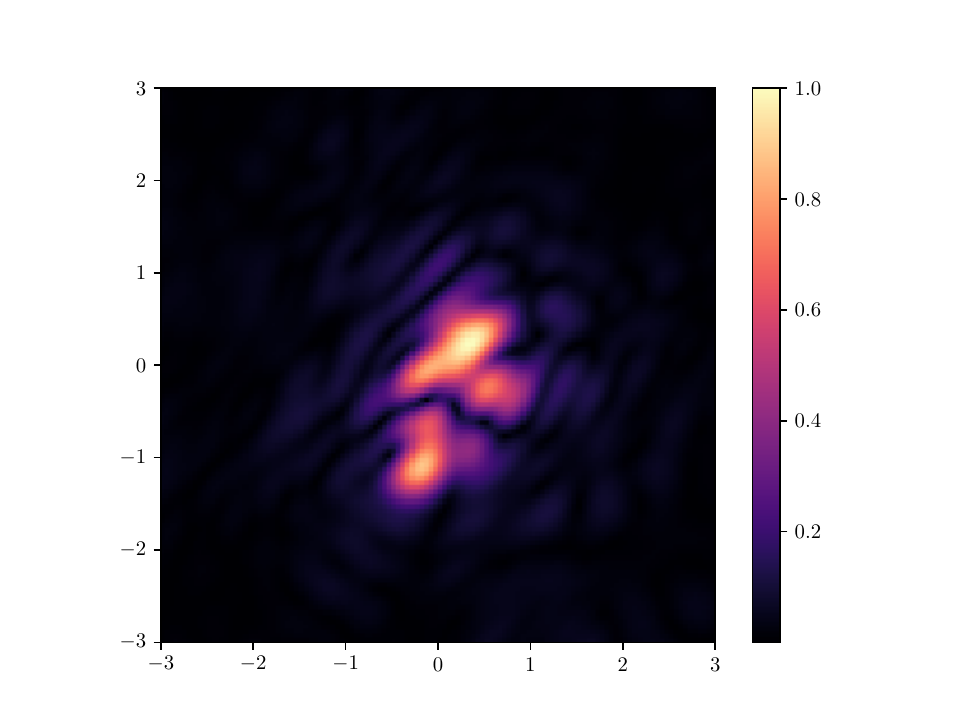}
		\caption{Normalised back-projection}
	\end{subfigure}~
	\begin{subfigure}[b]{0.45\textwidth}
		\centering
		\includegraphics[width=\columnwidth]{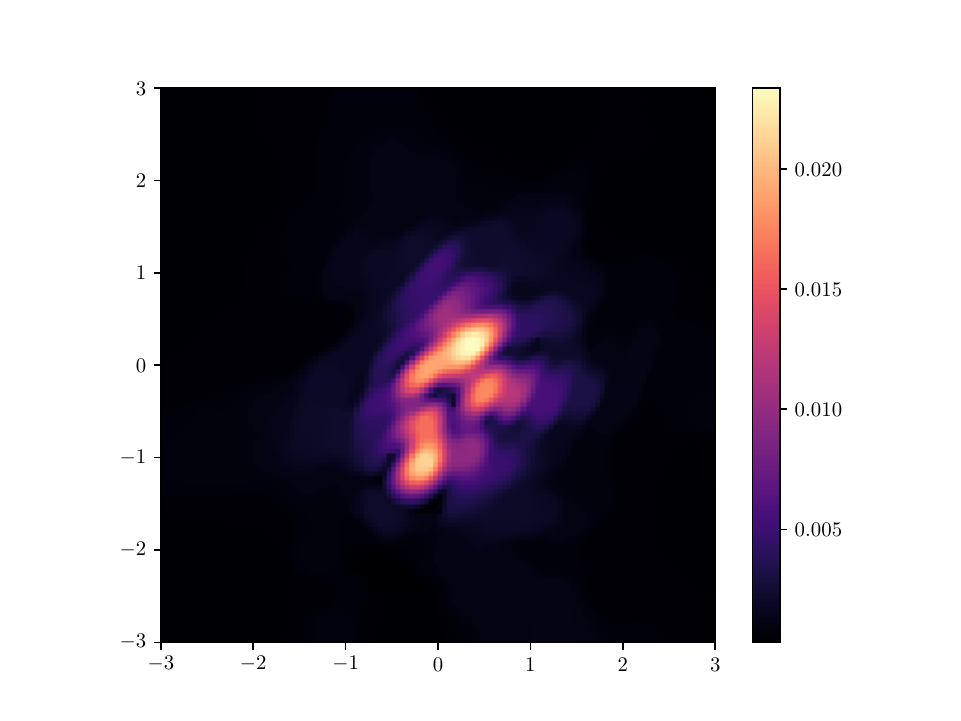}
		\caption{TV-regularised Born}
		\label{fig: recon2 tv born}
	\end{subfigure} \\
	\begin{subfigure}[b]{0.45\textwidth}
		\centering
		\includegraphics[width=\columnwidth]{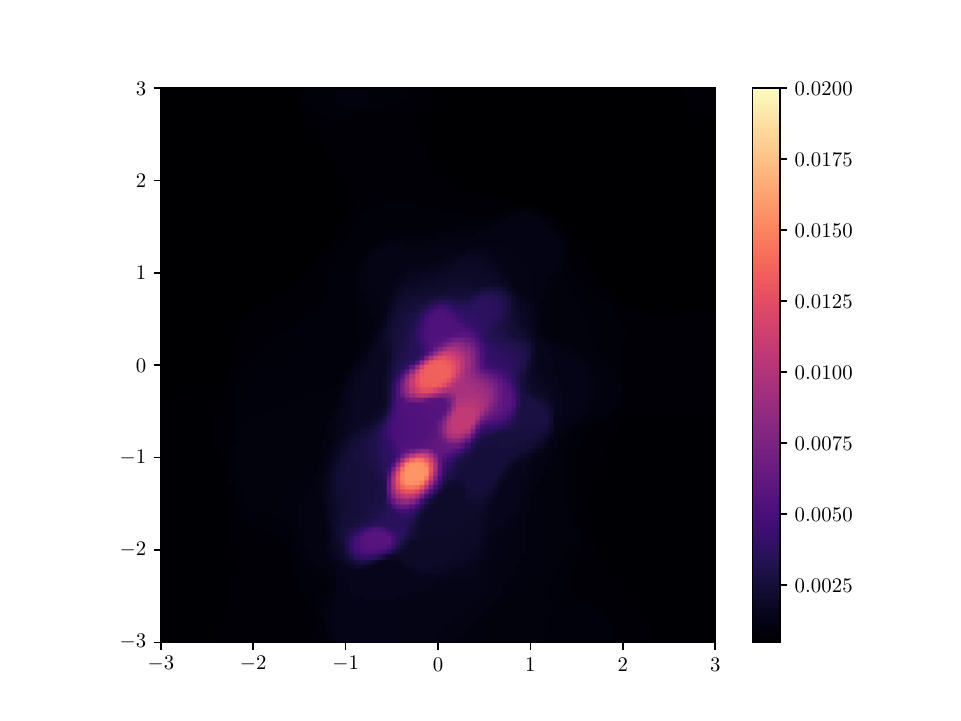}
		\caption{Iteration 1}
	\end{subfigure}~
	\begin{subfigure}[b]{0.45\textwidth}
		\centering
		\includegraphics[width=\columnwidth]{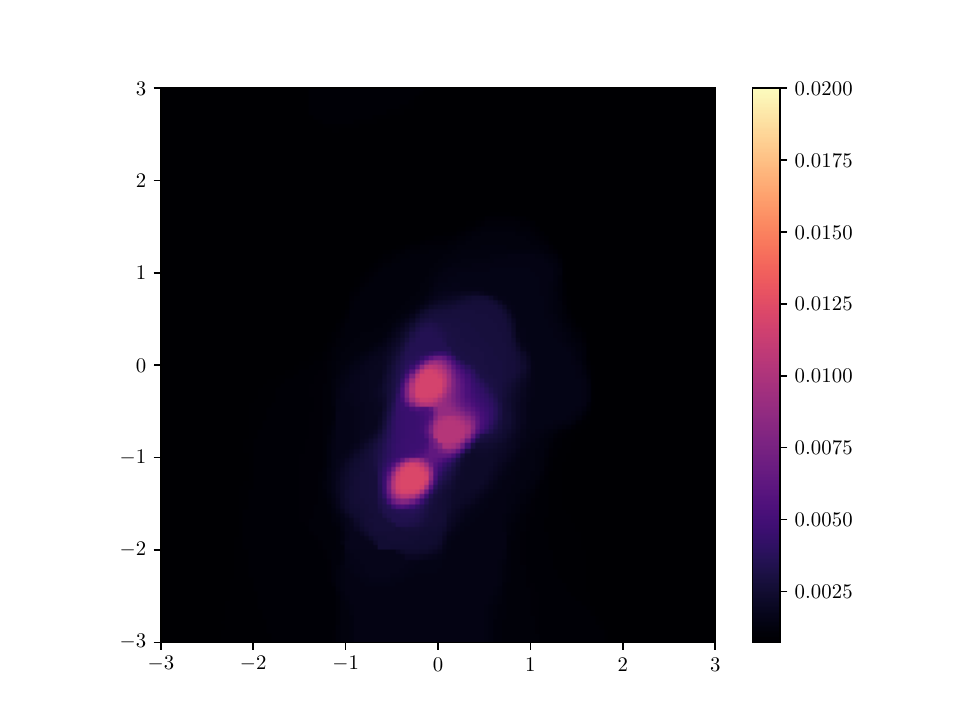}
		\caption{Iteration 5}
		\label{fig: recon2 final}
	\end{subfigure}
\caption{Magnitude of reflectivity for (a) normalised back-projection image, (b)  TV-regularised reconstruction with the standard Born approximation model, and the through-wall reconstruction \(|\bar{\lvec{v}}(\lele{m}^{[k]})|\) for (c) the 1\textsuperscript{st} and (d) the 5\textsuperscript{th} and final iterate \(\lele{m}^{[k]}\) during the BFGS optimisation of unknown permittivity \(\lele{m}=\epsilon_r\).  The final computed permittivity value is \(\epsilon_r=2.09\).}
\label{fig: 2d reconstruction unknown perm}
\end{figure}

\begin{figure}[ht]
	\centering
	\includegraphics[width=0.5\textwidth]{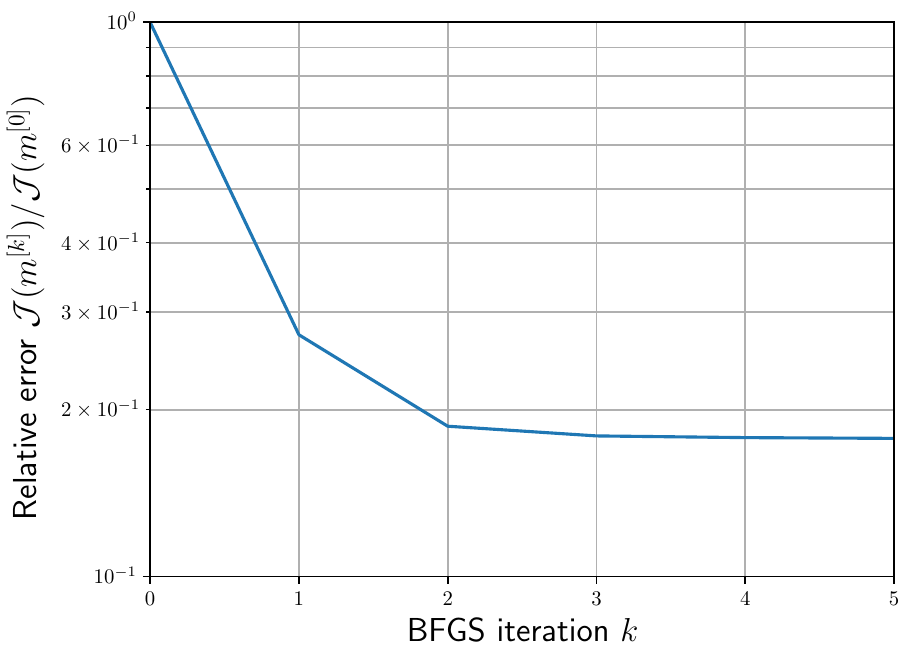}
	\caption{Convergence of BFGS optimisation of nuisance parameter \(\lele{m}=\epsilon_r\).}
	\label{fig: 2d recon convergence}
\end{figure}

The reconstruction results are shown in \figurename~\ref{fig: 2d reconstruction unknown perm}. It is clear that the targets behind the wall are rapidly focused in a small number of BFGS iterations when compared to both the standard SAR back-projection image and a TV reconstruction using the standard Born forward model (and otherwise the same parameters), again being much clearer that there are three distinct objects against a much cleaner background free of clutter.  We can see from \figurename~\ref{fig: 2d recon convergence} that the outer level of optimisation quickly stagnates after a small number of iterations.  The final calculated permittivity value is \(\epsilon_r=2.09\), which is reasonably accurate however, since this is a nuisance parameter, its accuracy matters less than achieving a well-focused image (i.e., the desire is not to know the specific material properties of a building, but what is inside it).

A lower peak value of the targets is observed in the reconstruction using the through-wall model versus the standard SAR model, though signal-to-background/clutter and target isolation is observed to be improved.  This is likely in part due to the choice of regularisation, with TV promoting step-like images allowing the target response to ``step up'' above the clutter in \figurename~\ref{fig: recon2 tv born}. There may also be some interference of out-of-phase neighbouring pixels (within the same effective resolution cell) in \figurename~\ref{fig: recon2 tv born} in trying to fit the data, due to the complex sidelobe and clutter pattern of the multi-static data, resulting in a greater overall (though misplaced and misformed) peak response.

\subsection{Reconstruction with an approximately known obstacle} \label{sec: recon 3}
\begin{figure}[h]
	\centering
	\includegraphics[width=0.75\textwidth, trim={0 2cm 0 2cm}, clip]{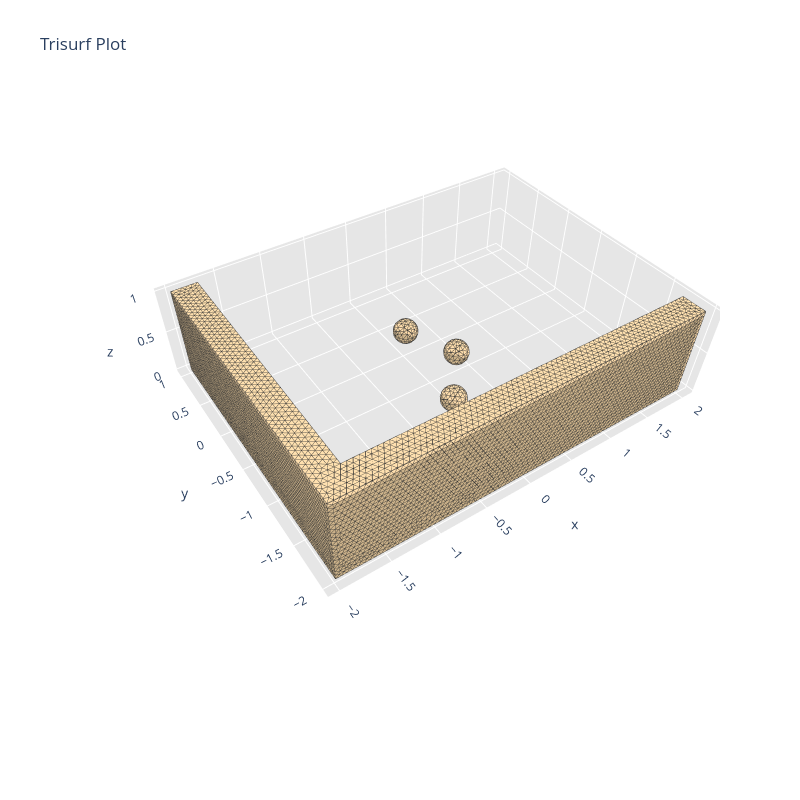}
	\caption{Data simulation domain.}
	\label{fig: true geom}
\end{figure}

Finally, we simulate data for three spheres of radius \qty{0.125}{\meter}, relative permittivity \(\epsilon_r=5.0\), and conductivity \qty{1.0}{\micro\siemens\per\meter}, behind a \qty{0.27}{\meter} thick non-conductive corner wall of permittivity \(\epsilon_r=2.85\), using BEMPP-CL\cite{betcke2021bempp} (i.e., we use a full wave model). The simulation domain is shown in \figurename~\ref{fig: true geom}.  As previously, \qty{20}{\percent} Gaussian noise is added to the data.

For the reconstruction we assume imperfect knowledge of the wall, assuming it is \qty{0.3}{\meter} thick, and with origin offset by \qty{-2}{\centi\meter} in both the \(x\)- and \(y\)-coordinates compared with the true obstacle.  We also use the full through-wall model for these reconstructions, i.e.
\begin{equation}
	\lvec{d}_{\mathrm{phd}} = \lmat{A}(\lvec{m})\lvec{v} + \mathcal{F}_0(\lvec{m}).
\end{equation}
We allow 10 iterations of BFGS to optimise for the permittivity of the assumed corner wall, which results in a value of \(\epsilon_r = 2.78\).  It makes sense that this is below the true permittivity, since the true wall is thinner than assumed, so the true value would result in too great an effective time delay.

\begin{figure}
	\centering
		\begin{subfigure}[b]{0.45\textwidth}
				\centering
				\includegraphics[width=\columnwidth]{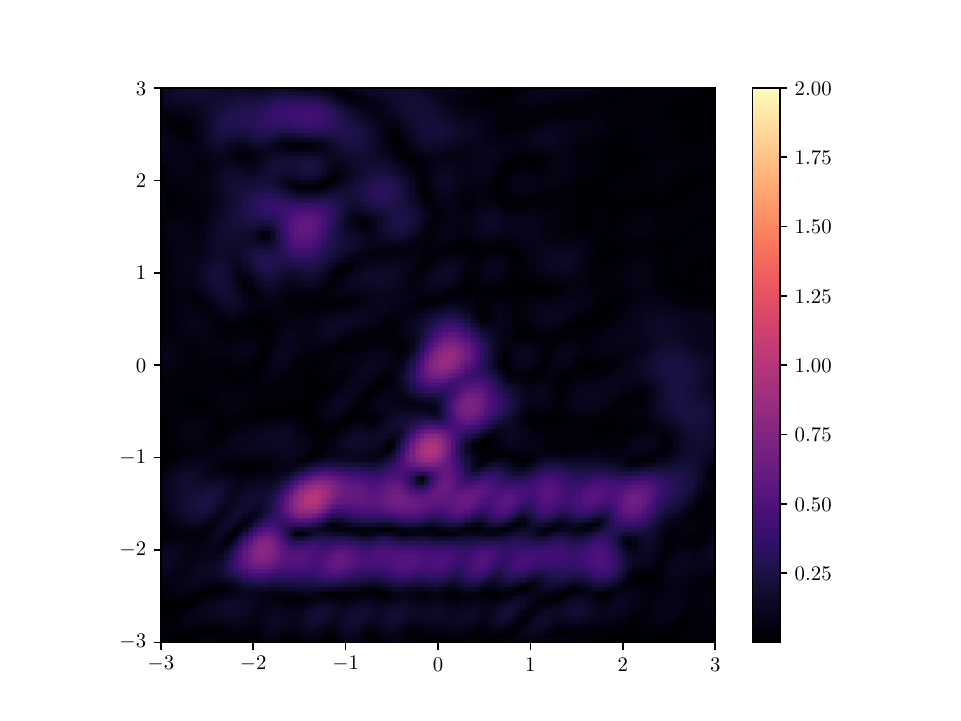}
				\caption{Standard SAR back-projection}
				\label{fig: std bpim}
			\end{subfigure}~
		\begin{subfigure}[b]{0.45\textwidth}
				\centering
				\includegraphics[width=\columnwidth]{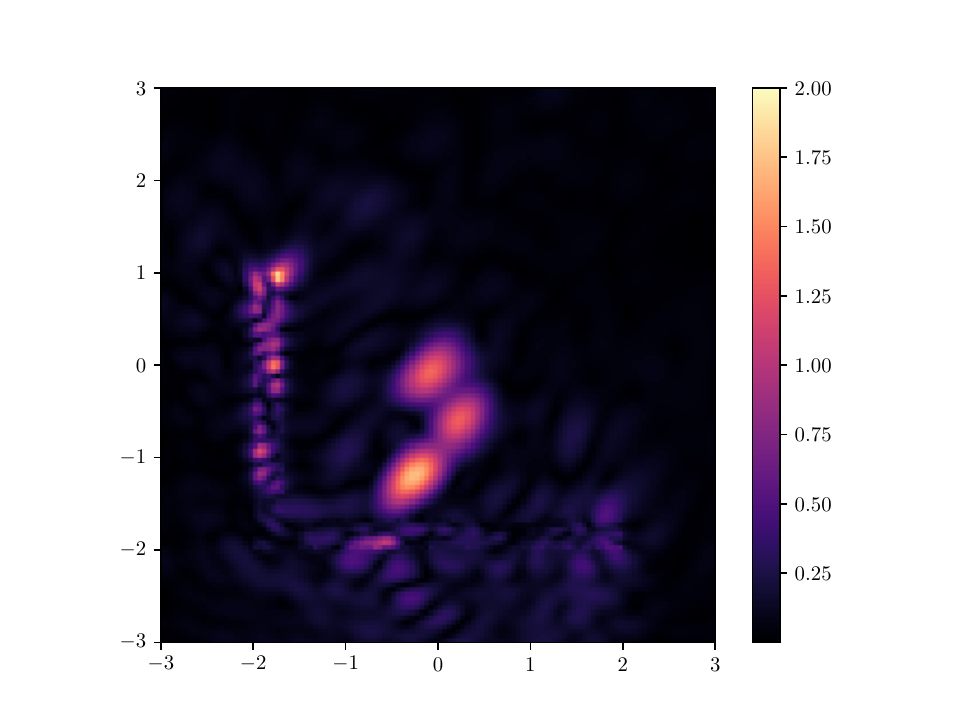}
				\caption{Through-wall Back-projection}
				\label{fig: tw bpim}
			\end{subfigure}~\\
		\begin{subfigure}[b]{0.45\textwidth}
			\centering
			\includegraphics[width=\columnwidth]{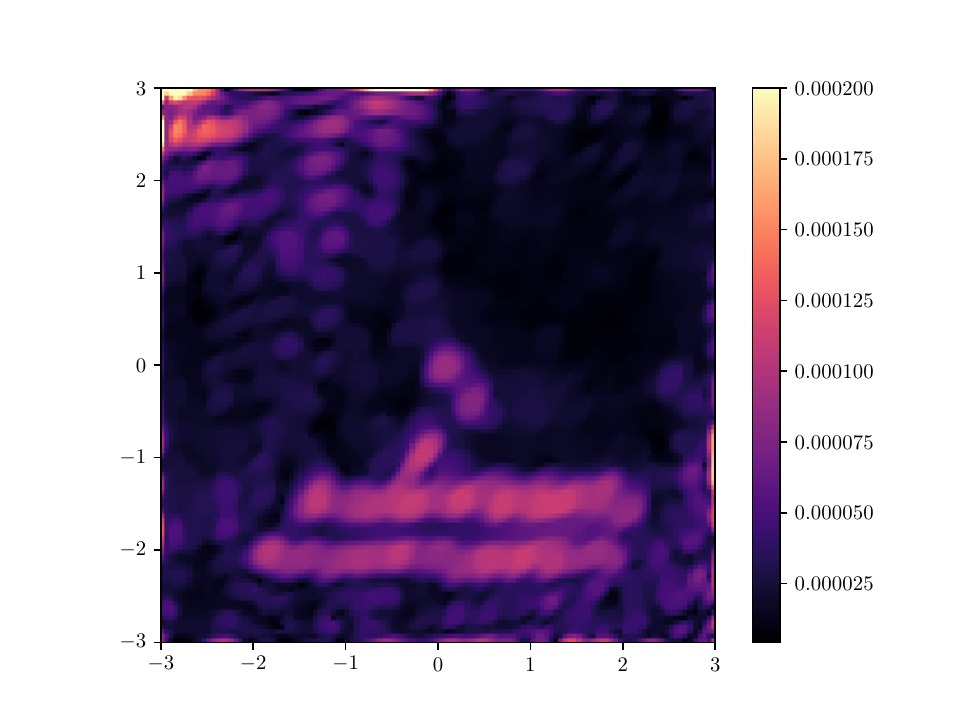}
			\caption{Standard SAR TV-regularised}
			\label{fig: std TV}
		\end{subfigure}~
		\begin{subfigure}[b]{0.45\textwidth}
			\centering
			\includegraphics[width=\columnwidth]{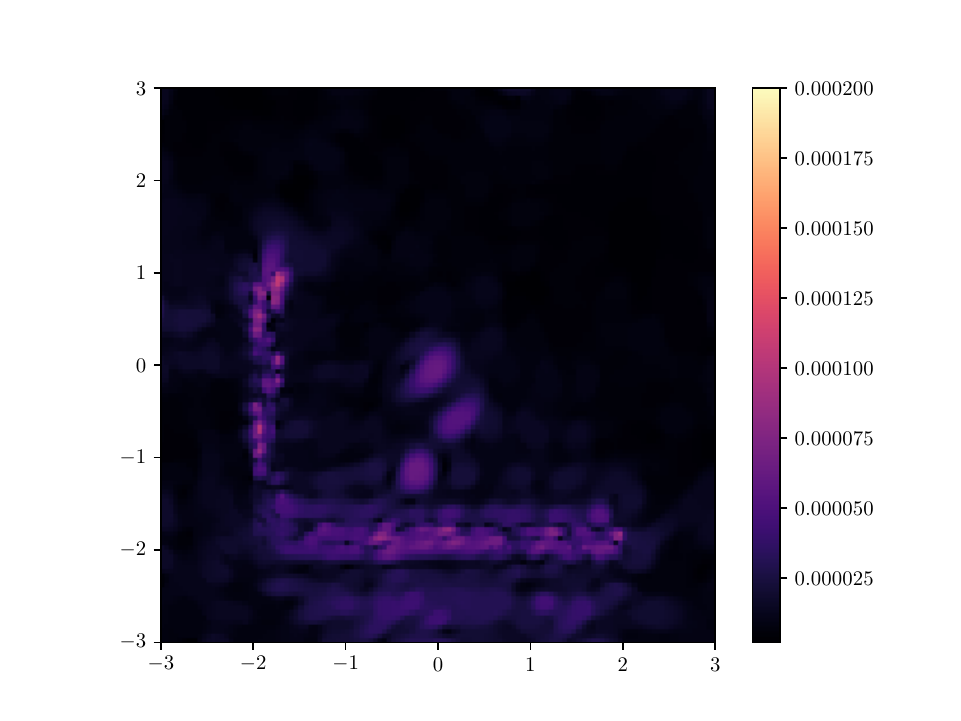}
			\caption{Through-wall TV-regularised}
			\label{fig: recon3 TV}
		\end{subfigure}~\\
		\begin{subfigure}[b]{0.45\textwidth}
			\centering
			\includegraphics[width=\columnwidth]{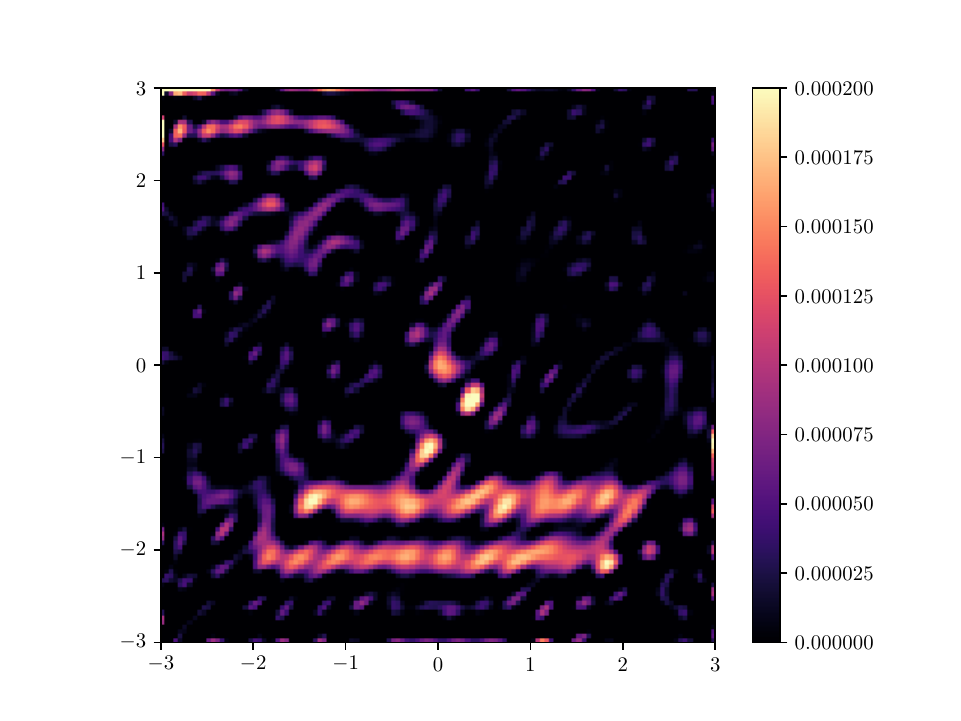}
			\caption{Standard SAR L1-regularised}
			\label{fig: std L1}
		\end{subfigure}~
		\begin{subfigure}[b]{0.45\textwidth}
			\centering
			\includegraphics[width=\columnwidth]{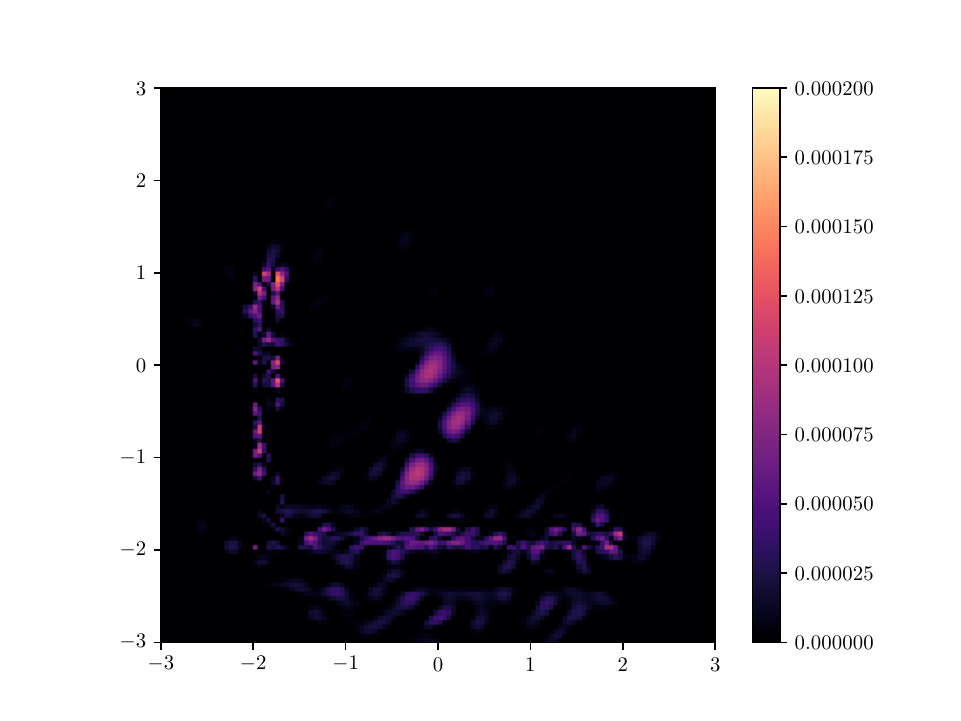}
			\caption{Through-wall SAR L1-regularised}
			\label{fig: recon3 L1}
		\end{subfigure}
		\caption{Magnitude of reflectivity images for (a) back-projection with the standard SAR model, (b) back-projection with the through-wall model, (c) TV-regularised reconstruction with the standard SAR model, (d) TV-regularised reconstruction with the through-wall model, (e) L1-regularised reconstruction with the standard SAR model, and (e) L1-regularised reconstruction with the through-wall model.  (a) and (b) are both normalised to the maximum value of the standard SAR back-projection.  (b), (d) and (e) are all obtained using the optimised permittivity value \(\epsilon_r=2.78\). }
		\label{fig: recon3 set}
\end{figure}

The reconstruction results are shown in \figurename~\ref{fig: recon3 set}. The left column shows the back-projection, TV-regularised reconstruction, and L1-regularised reconstruction results using the standard SAR model, and the right column the same using the through-wall model.  \figurename~\ref{fig: recon3 TV} is the implicitly defined reflectivity resulting from the optimisation of permittivity, with \figurename~\ref{fig: tw bpim}~and~\ref{fig: recon3 L1} being produced directly using the already resolved permittivity value \(\epsilon_r=2.78\) without further BFGS optimisation.
  
Immediate observations are the effects of interactions with the wall itself when using the standard SAR model in both back-projection and reconstructions: both overlay of two strong reflections of the wall aligned to the \(x\)-direction, and the effect of multiple reflections off of the wall itself projecting into the top-left of these images.  The latter are particularly poorly dealt with in the regularised least-squares reconstructions, most likely since these multiple reflections are not in the range of the forward operator.  The results may well be improved by first filtering data (likely challenging with imperfect knowledge of the structure and scene), but we find these artefacts are well resolved by the output of our reconstruction method shown in \figurename~\ref{fig: recon3 TV}. This demonstrates some viability in our method without the very strong prior of a perfectly known obstacle geometry.  

Total Variation appears to be a poor choice in reducing artefacts for the standard SAR model reconstruction shown in \figurename~\ref{fig: std TV} (it rather appears only to join these artefacts together), so we also compare L1-regularised reconstructions in Figures~\ref{fig: std L1}~and~\ref{fig: recon3 L1}, i.e. with penalty \(\mathcal{R}(\lvec{z})=\|\lvec{z}\|_1\).   This again does relatively little to suppress artefacts using the standard SAR forward model, but does result in a very low level of artefacts when using the through-wall model and arguably the most visually appealing reconstruction of the four.  This highlights the need to tailor the particular regularisation used to the particular real-world scenario and radar system used for through-wall imaging. TV alone may prove to be less useful in more realistic scenarios containing complex objects, but some combination of L1 plus TV may be beneficial.

Finally, we also show the back-projection images for the standard SAR model in \figurename~\ref{fig: std bpim}, and for the through-wall model (again with \(\epsilon_r=2.78\)) in \figurename~\ref{fig: tw bpim}.  Both are normalised to the maximum value of the standard SAR back-projection. The brighter scatterers in \figurename~\ref{fig: tw bpim} demonstrates the increased signal-to-clutter observed.  The quality of backprojection image \figurename~\ref{fig: tw bpim} highlights that the through-wall model proposed may also be used in producing standard SAR imagery, potentially forgoing costly reconstruction schemes, provided one has a good enough estimate of the obstacle's properties. For more complex scenes some level of optimisation or manual testing of parameters may be required, though this is to be determined using real data.

\FloatBarrier
\section{Discussion and Conclusions}
We have presented a partially linearised reconstruction method for through-wall SAR image formation, able to resolve through-wall multiple scattering effects. This uses a more representative model of the full physical scattering process than has previously been applied to SAR, whilst being significantly less computationally expensive than the full-wave type reconstruction formulations which have been applied to related through-wall and sub-surface radar imaging problems with different geometries and considerations. Thus, it provides a method which may be practical for real-world applications of stand-off imaging now. That is, through-wall transmission effects (including multiple scattering through walls) are included for an arbitrary building, walls, or other obstacle.
 The method is also applicable to quite general data collections geometries.

The computational efficiency is gained by simulating only the through-wall effects with a full-wave model using a boundary element method, which provides numerical approximations to Green's functions for the scene.  A strong prior -- that we are imaging through a wall with (approximately) known geometry -- allows us to reduces this simulation problem to one of primarily freespace transmission.  Coupling the BEM simulation with a ROM provides near real-time formation of the linear SAR simulation model (Born approximation) for given wall properties.  These features are desirable in security applications, where an image may need to be formed quickly, and there may be limitations in the collection geometries which can be achieved from longer ranges, but it may still be practical to commit the computational resources required in the \textit{offline} stage during any mission planning, saving time later.

The method has been tested in the case of an approximately known obstacle with unknown permittivity.  More unknowns about the structure (such as geometrical properties) may easily be incorporated in the same framework to be resoled as part of the imaging process, at increased computational and memory cost of the \textit{offline} stage, though we have not demonstrated this capability here.  Further research is required to understand the number and type of unknowns which can be resolved, and what form of ``rich'' multi-dimensional data (here using multi-statics) is required to do so.  Moreover, we have demonstrated the method only for 2D image formation. However, since the model itself is 3D, this means any out-of-plane through-wall transmission effects are resolved even in 2D images, and it may be directly applied to 3D datasets to also resolve objects in height.

While the model does capture more of the physics, there are still several features which it currently does not.  First and foremost, we have still made the common assumption of scalar waves -- clearly this is incorrect, especially where there is multiple scattering, but it is a standard approximation to make in SAR.  This is easily updated  by replacing the BEMPP Helmholtz solver with the equivalent for Maxwell's equations -- which is essentially a black-box component for the reconstruction method -- albeit with increased computational cost.  One would also have to consider more carefully how polarization data should be handled. Would it be more useful to form multiple images in different polarizations to allow analysts to interpret (as is often the case) -- and should these images be separated by the (multi-static) antenna polarization data channels, or by the polarization of the incident through the wall.  

Alternatively, do we include polarization changing effects in the Born scattering model, resulting in a tensor-valued (``rich'') reconstruction problem in which we would have to make some assumptions about the effective shapes we could represent each voxel having (i.e., different polarization changes would be due to the extent to which the scatterers are \emph{not} point-like).  The latter may be more appropriate, and would be captured by default in a full-wave reconstruction method, but the result could have the unwanted effect of making the information harder to interpret by the user.  An appropriate model is therefore to be determined.

Relatedly, we have also not included a floor in the simulations provided here.  Double bounces off of the floor and a wall may be prominent features in the data, forming an edge reflector which will direct significant energy back to a mono-static radar when near oblique to the wall, or the equivalent bi-static geometries.  This feature is perhaps most readily introduced by using an analytic Green's function for a half-space (or 2-layered space) in the boundary integral formulation, which would avoid needing to mesh additional features at the expense of a more involved evaluation of the equivalent boundary element integrals to (\ref{eq: discrete operators}).  Otherwise, additional structural elements (more walls, roofs, windows) are straightforward to include in our formulation, potentially adding additional nuisance parameters if their structure is not known exactly and increasing the computational cost of the offline stage.

Finally, we remark that further attention should be given to the optimisation schemes used in refining a practical fieldable method.  For the outer level of optimisation, understanding the degree to which convexity here might rely on particular multi-static geometries.  Given the small number of parameters, a derivative free or Bayesian approach to sampling these parameters may be more appropriate.  For the inner problem, the choice of regularisation would need to be further tailored to specific applications and radar systems, and more efficient optimisation algorithms such as PDHG or ADMM may be applied.  Alternatively, it may often prove to be most beneficial to simply use the through-wall forward model in producing back-projection images with an estimate of the scene properties, since this has very little computational overhead in the \textit{online} stage.

\ack
The authors would like to thank the Isaac Newton Institute for Mathematical Sciences, Cambridge, for support and hospitality during the programme ``Rich and Nonlinear Tomography: A Multidisciplinary Approach'', where work on this paper was undertaken. This work was supported by EPSRC grant no EP/R014604/1 and EP/V007742/1, and has made use of computational support by CoSeC, the Computational Science Centre for Research Communities, through CCPi.

Watson was supported by the Royal Academy of Engineering and the Office of the Chief Science Adviser for National Security under the UK Intelligence Community Postdoctoral Research Fellowship programme.

\begin{appendices}
	\numberwithin{equation}{section}
	\section{Boundary element formulation of the through-wall model}\label{ap: BEM}
	Here we provide the derivation of the boundary element formulation of the representation formula of the scattering problem (\ref{eq: Helmholtz inclusion}), namely
	\begin{equation}
		U(\mathbf{x},\omega) = \int_{\partial D} G_0(\mathbf{x},\mathbf{y},\omega)\psi(\mathbf{y}) \rmd\mathbf{y} - \int_{\partial D} \frac{\partial G_0(\mathbf{x},\mathbf{y},\omega)}{\partial \mathbf{n}}\phi(\mathbf{y}) \rmd\mathbf{y},
		\label{eq: ap representation formula}
	\end{equation}
	with jumps of the solution across the interface
	\begin{align}
		\psi &= \gamma_1^- u - \gamma_1^+ u, \\
		\phi &= \gamma_0^- u - \gamma_0^+ u,
	\end{align}
	following Haqshenas \textit{et al}\cite{Haqshenas21}.
	
	\subsection{Surface integral equations}
	By taking the Cauchy trace for the exterior \(\gamma^+\)/interior \(\gamma^-\), \(\gamma^{\pm}:=\left[\gamma_0^{\pm},\,  \gamma_1^{\pm}\right]^T\), of the representation formula (\ref{eq: ap representation formula}) it can be shown\cite{Betcke17, Betcke22, Haqshenas21, Claeys13} that in the exterior of \(D\) the wavefield satisfies
	\begin{equation}
		\left(\frac{1}{2}Id - \mathcal{A}\right)\gamma^+u= \gamma^+u^{sc},
		\label{eq: exterior bie}
	\end{equation}
	and in the interior of \(D\)
	\begin{equation}
		\left(\frac{1}{2}Id + \mathcal{A}\right)\gamma^-u = \gamma^- u,
		\label{eq: interior bie}
	\end{equation}
	Here, \(Id\) is the identity operator and \(\mathcal{A}\) is the \textit{Calder{\'o}n} operator
	\begin{equation}
		\mathcal{A} = \begin{bmatrix} -\mathcal{K} & \mathcal{V} \\ \mathcal{W} & \mathcal{K}'
		\end{bmatrix},
	\end{equation}
	with single-, double-, hypersingular- and adjoint double-layer boundary integral operators \(\mathcal{K}\), \(\mathcal{V}\), \(\mathcal{W}\) and \(\mathcal{K}'\) given by
	\begin{align}
		(\mathcal{V}\psi)(\mathbf{x})&:= \int_{\partial D} G_0(\mathbf{x},\mathbf{y},\omega)\psi(\mathbf{y}) \rmd\mathbf{y}, \quad  \\
		(\mathcal{K}\phi)(\mathbf{x}) & :=\int_{\partial D} \frac{\partial G_0(\mathbf{x},\mathbf{y},\omega)}{\partial \mathbf{n}(\mathbf{y})}\phi(\mathbf{y}) \rmd\mathbf{y}, \\
		(\mathcal{W}\phi)(\mathbf{x}) &:= -\frac{\partial}{\partial\mathbf{n}(\mathbf{x})}\int_{\partial D}\frac{\partial G_0(\mathbf{x},\mathbf{y},\omega)}{\partial\mathbf{n}(\mathbf{y})}\psi(\mathbf{y})\rmd\mathbf{y}, \\
		(\mathcal{K}'\psi)(\mathbf{x})&: = \frac{\partial}{\partial \mathbf{n}(\mathbf{x})}\int_{\partial D}  G_0(\mathbf{x},\mathbf{y},\omega)\phi(\mathbf{y}) \rmd\mathbf{y},
	\end{align}
	for \(\mathbf{x}\in\partial D\) in each case.  Combining (\ref{eq: interior bie}) with (\ref{eq: exterior bie}) we have that\cite{Haqshenas21,Claeys13}
	\begin{equation}
		\left(\mathcal{A}_{(k_0)} + \mathcal{A}_{(k_{D})}\right)\gamma^+u = \gamma^+u^{in},
		\label{eq: combined bie}
	\end{equation}
	which completes the reformulation of the exterior model (\ref{eq: Helmholtz inclusion}) in an infinite domain with a boundary integral equation over the surface of the inclusion \(\partial D\).  The solution to (\ref{eq: combined bie}) provides the Dirichlet and Neumann traces of the scattered wavefield on the surface \(\partial D\), with which we can evaluate \(u\) at any point via the representation formula (\ref{eq: representation formula}).

	\subsection{Finite element formulation}
	The boundary integral formulation (\ref{eq: combined bie}) may be discretised with the Galerkin method.  This first involves forming a triangular mesh \(\partial D_h\) of the boundary \(\partial D\) with \(n\) nodes \(\mathbf{x}_j\), \(j=1,2,\ldots n\), in this case piecewise linear elements \(\varphi_j\), which at the nodes take values
	\begin{equation*}
		\varphi_j(\mathbf{x}_i) = \left\{\begin{array}{ll} 1, & i=j \\ 0, & i\neq j.\end{array}\right.
	\end{equation*}
	Using the elements \(\varphi_j\) to approximate the solutions \(\gamma\), discretising the weak formulation of (\ref{eq: combined bie}) provides us with the Boundary Element (BEM) linear system
	\begin{equation}
		\lmat{A}\lvec{u} = \lvec{b},
		\label{eq: bempp system}
	\end{equation}
	where \([\lvec{A}]_{ij} = \left<(\mathcal{A}_{(k_0)} + \mathcal{A}_{(k_{D})})\varphi_j, \varphi_i\right>\), \(b_i=\left<\gamma^+u^{in},\varphi_i\right>\).  Computing these matrix entries involves evaluating
	\begin{subequations}
		\begin{align}
			[\lmat{V}]_{ij}&:= \int_{\partial D} \varphi_i(\mathbf{x})\int_{\partial D} G_0(\mathbf{x},\mathbf{y},\omega)\varphi_j(\mathbf{y}) \rmd\mathbf{y}\rmd\mathbf{x} \\
			[\lmat{K}]_{ij} & :=\int_{\partial D}\varphi_i(\mathbf{x})\int_{\partial D} \frac{\partial G_0(\mathbf{x},\mathbf{y},\omega)}{\partial \mathbf{n}(\mathbf{y})}\varphi_j(\mathbf{y}) \rmd\mathbf{y}\rmd\mathbf{x}, \\
			[\lmat{W}]_{ij} &:= -\int_{\partial D}\varphi_i(\mathbf{x})\frac{\partial}{\partial\mathbf{n}(\mathbf{x})}\int_{\partial D}\frac{\partial G_0(\mathbf{x},\mathbf{y},\omega)}{\partial\mathbf{n}(\mathbf{y})}\varphi_j(\mathbf{y})\rmd\mathbf{y}. \label{eq: discrete hypersingular}
		\end{align}
		\label{eq: discrete operators}
	\end{subequations}
	Note that the discretisation of the adjoint double-layer operator \(\mathcal{K}'\) is given simply by the adjoint of \(\lmat{K}\), and that integration by parts of (\ref{eq: discrete hypersingular}) yields a weakly singular integral.  Evaluation of these integrals and assembly of the linear system and source terms in (\ref{eq: bempp system}) are handled fully by the BEMPP-CL library \cite{betcke2021bempp}, which includes suitable quadrature rules for the singular integrals as well as sharp-edged meshes.  The resulting systems are often solved via a pre-conditioned iterative method such as GMRES.  However, since we need to solve many such systems at each given frequency for different source locations (along synthetic apertures) we instead use a direct solver, calculating and storing an LU factorisation of the dense system matrices for reuse with each subsequent source term, which has been feasible for the size of problems considered so far.

	\section{Reconstruction details}\label{ap: reconstruction}
	\subsection{Objective function and derivatives} 
	\label{ap: reconstruction derivs}
		In this work, we use a variable projection formulation of the objective function in reflectivity \(\lvec{v}\) and wall parameters \(\lvec{m}\),
		\begin{subequations}
			\begin{align}
				\tilde{\lvec{m}} &= \argmin_{\lvec{m}} \bar{\mathcal{J}}(\mathbf{m}) \label{eq: appendix projected base}\\
				\bar{\mathcal{J}}(\lvec{m})&:= \mathcal{J}(\lvec{m};\bar{\lvec{v}}(\lvec{m})) \label{eq: appendix projected m}\\
				\bar{\lvec{v}}(\lvec{m}) &:= \argmin_{\lvec{v}} \mathcal{J}(\lvec{m},\lvec{v}). \label{eq: appendix projected v}
			\end{align}
			\label{eq: appendix reconstruction problem projected}
		\end{subequations}		
		Derivatives of (\ref{eq: appendix projected m}) with respect to \(\lvec{m}\) are given by \cite{aravkin2012estimating}
		\begin{equation}
			\nabla \bar{\mathcal{J}}(\lvec{m}) = \left.\nabla_m \mathcal{J}(\lvec{m},\lvec{v})\right|_{\lvec{v}=\bar{\lvec{v}}(\lvec{m})},
		\end{equation}
		and
		\begin{equation}
			\nabla^2 \bar{\mathcal{J}} = \Bigl(\nabla^2_{\lvec{m}\lvec{m}}\mathcal{J}(\lvec{m},\lvec{v}) - \nabla^2_{\lvec{m}\lvec{v}}\mathcal{J}(\lvec{m},\lvec{v})(\nabla_{\lvec{vv}}\mathcal{J}(\lvec{m},\lvec{v}))^{-1} \nabla^2_{\lvec{mv}}\mathcal{J}(\lvec{m},\lvec{v}) \Bigr|_{\lvec{v}=\bar{\lvec{v}}(\lvec{m})},
			\label{eq: implicit Hess}
		\end{equation}
		with terms involving \(\nabla_{\lvec{v}}\mathcal{J}(\lvec{m},\lvec{v})\bigr|_{\lvec{v}=\bar{\lvec{v}}(\lvec{m})}\) having dropped out since \(\bar{\lvec{v}}\) is (possibly approximately) a minimum.  It has been noted in the related problem of differential semblance that where the implicit problem (\ref{eq: reconstruction problem projected}) is not solved exactly that these terms may in fact be far from zero\cite{symes1991velocity}, and a correction term may be warranted.  Since such correction terms can be expensive to calculate, for this work, we will assume that (\ref{eq: implicit v}) has been solved sufficiently accurately.
		
		Evaluating the gradient terms, we have
		\begin{align}
			\frac{\partial \mathcal{J}}{\partial m_j} &= \left(\frac{\partial \mathcal{F}_0}{\partial m_j} + \frac{\partial \mathcal{F}_1}{\partial m_j}, \mathcal{F}_0 +\mathcal{F}_1- \lvec{d}\right) \\
			\left[\frac{\partial \mathcal{F}_0}{\partial m_j}\right]_i &= a(\omega_i)\rme^{-\rmi\omega_i R_{0,i}/c} \left[H_{P0}^M\Sigma_{P0}^M\frac{\partial \lvec{g}_{P0}}{\partial m_j}\right]_i,
		\end{align}
		and for field scattered by objects in the scene
		\begin{align}
			\lmat{J}_i:=\frac{\partial \mathcal{F}_1}{\partial m_j} =& \frac{\partial A}{\partial m_j} \lvec{v}, \qquad \lmat{J}:=[\lmat{J}_1,\ldots,\lmat{J}_{N_m}] \label{eq: Jacobian} \\
			\left[\frac{\partial \mathcal{F}_1}{\partial m_j}\right]_i=& a(\omega_i)\rme^{-\rmi\omega_i R_{0,i}/c}\cdot  \nonumber \\
			&\sum_j\Biggl\{\Biggl( \left(G_0(\mathbf{x}_j, \mathbf{y}_{T,i}) + \left[H_{P1}^M\Sigma_{P1}^M \lvec{g}_{P1}(\omega_i, s_i)\right]_{Ri,j} \right) \nonumber \\
			&\qquad\qquad\left[H_{P1}^M\Sigma_{P1}^M \frac{\partial \lvec{g}_{P1}}{\partial m_i}(s_i,\omega_i)\right]_{Ti,j}   \nonumber\\
			& \qquad\quad +  \left(G_0(\mathbf{x}_j, \mathbf{y}_{R,i}) + \left[H_{P1}^M\Sigma_{P1}^M \lvec{g}_{P1}(\omega_i, s_i)\right]_{Ti,j} \right) \nonumber \\
			&\qquad\qquad\left.\left[H_{P1}^M\Sigma_{P1}^M \frac{\partial \lvec{g}_{P1}}{\partial m_i}(s_i,\omega_i)\right]_{Ri,j}\right) v_j \Biggr\}.
		\end{align}
		The gradients \(\nabla_{\lvec{m}}\lvec{g}_P\) are available directly as the gradients of the interpolation splines in an appropriate interpolation library. 
		
		From (\ref{eq: implicit Hess}) and (\ref{eq: Jacobian}), we can see the Hessian of \(\bar{\mathcal{J}}\) will involve matrix products of both first and second derivatives of \(\lmat{A}\) with respect to \(\lmat{m}\), as well as a term \(\left(\lmat{A}^*\lmat{A}\right)^{-1}\).  Since we have only a small number of nuisance parameters \(\lmat{m}\), but the reflectivity scene \(\lmat{v}\) is potentially quite large, this will be expensive to calculate directly relative to the scale of the problem.  As such we have used a Hessian-free method to recover these parameters.
		
		\subsection{Algorithmic details for the inner problem in reflectivity}\label{ap: algorithm reflectivity}
		The inner optimisation problem in reflectivity \(\bar{\lvec{v}}(\lvec{m})\) has been regularised by Total Variation applied to the reflectivity function, and solved via the FISTA algorithm. Isotropic TV is defined by
		\begin{equation}
			TV(\lvec{v}) := \|D \lvec{v}\|_{2,1} = \sum_{ij} \sqrt{(\lmat{D}_x \lvec{v})^2_{ij} + (\lmat{D}_y \lvec{v})^2_{ij}}
		\end{equation}
		in two dimensions, with \(D\) the discrete Laplacian operator \(\nabla\) and \(D_{p}\) the discrete derivative operator in the \(p\)\textsuperscript{th} coordinate.  The proximal map for \(TV(|\cdot|))\) can be shown to be given by\cite{Guven16}
		\begin{equation}
			\prox_{TV(|\cdot|)}(\lvec{y}) = \prox_{TV}(\lvec{r})\circ\exp^{\rmi\angle(\lvec{y})},
			\label{eq: prox TV}
		\end{equation}
		where \(\lvec{r}=|\lvec{y}|\) and \(\circ\) is the element-wise vector product, and the right-hand-side of (\ref{eq: prox TV}) can be computed efficiently via Fast Gradient Projection\cite{Beck09} (or FISTA). For the results in Section~\ref{sec: recon 3} we also apply L1 regularisation, i.e. \(\mathcal{R}(\lvec{z}) = \|\lvec{z}\|_1\). This has proximal map with i\textsuperscript{th} component given by
		\begin{equation}
			\prox_{\lambda\|\cdot\|}(\lvec{y})_i = \max(|\lele{y}_i|-\lambda,0)\exp^{\rmi\angle \lele{y}_i}.
			\label{eq: prox L1}
		\end{equation}
		
		Equipped with these proximal maps, several proximal splitting methods exist for the efficient solution of (\ref{eq: implicit v}). We implement (\ref{eq: implicit v}) in the CCPi Core Imaging Library\cite{jorgensen2021core, papoutsellis2021core}, and for simplicity use the implementation of FISTA therein.
		
		In solving the inner problem \(\bar{\lvec{v}}(\lvec{m})\), FISTA is warm started with the previously found solution \(\bar{\lvec{v}}'(\lvec{m}')\), reducing the number of iterations required for small changes \(\lvec{m}-\lvec{m}'\).  To realise this reduction, the algorithm is halted when either the relative changes in objective value \(\mathcal{J}^{[k]}\) or iterate \(\lvec{v}^{[k]}\) are less than some tolerances,
		\begin{subequations}
			\begin{align}
				\frac{\mathcal{J}(\lvec{m},\lvec{v}^{[k-1]})-\mathcal{J}(\lvec{m},\lvec{v}^{[k]})}{\mathcal{J}(\lvec{m},\lvec{v}^{[k-1]})}<r_{tol}, \\
				\frac{\|\lvec{v}^{[k-1]}-\lvec{v}^{[k]}\|_2^2}{\|\lvec{v}^{[k-1]}\|_2^2}<v_{tol}.
			\end{align}
		\end{subequations}
		Otherwise, FISTA is allowed to continue for a predetermined maximum number of iterations.
		
		\subsection{Algorithm details for the outer problem in nuisance parameters}\label{ap: algorithm outer}
		The outer optimisation problem in \(\lvec{m}\) may be solved via any suitable non-linear black-box optimisation scheme.  To avoid direct computation of Hessian matrices, we have chosen the BFGS scheme\cite{nocedal1999numerical}\footnote[1]{For convenience we actually use an existing l-BFGS implementation, but by setting the memory parameter equal to the total number of iterations and maintaining the same initial Hessian approximation the two are equivalent, albeit with some very minor increased computational cost through very small matrix vector products.}.  Since we have only a small number of nuisance parameters  \(\lvec{m}\), we calculate the initial approximate Hessian inverse as \(\lmat{H}_0 = \left(J_0^*J_0\right)^{-1}\), where \(J_0\) is the Jacobian of \(\mathcal{F}\) with respect to \(\lvec{m}\) evaluated at the initial estimate \(\lvec{m}_0\).
		
		We have not found it necessary to implement bound constraints such as relative permittivity \(\epsilon_r\geq1\) as our experience is that \(\bar{\mathcal{J}}(\lvec{m})\) will begin to increase greatly away from the true parameters \(\lvec{m}_{true}\).  Indeed, extrapolation of the POD outside the region with snapshots may add such numerical error as to act as an additional penalty. This should however be given proper consideration for future work, in particular if the method is adapted for a more general parameterisation of the wall/obstacle.  
		
		Equally, in these initial numerical experiments we have not found any need for additional regularisation of \(\lvec{m}\).  In particular, since \(\lvec{m}\) is a nuisance parameter, the motivating purpose is not to determine the wall structure to a particular accuracy (as it might be in e.g. a non-destructive testing application), but to resolve it \textit{accurately enough} to form a well focused image \(\lvec{v}\). Moreover, the setup already enforces very strong prior knowledge that the wall is indeed ``wall shaped''.  Again, this may need further consideration for a more general parameterisation.
		
\end{appendices}

\section*{Data availability statement}
The data that support the findings of this study are available upon request from the authors.

\FloatBarrier

\printbibliography

\end{document}